\pdfoutput=1
\documentclass[reqno]{amsart}
\usepackage{algorithm}
\usepackage{amsmath, amssymb}
\usepackage[foot]{amsaddr}
\usepackage{bm}
\usepackage{color}
\usepackage{csquotes}
\usepackage{gensymb}
\usepackage[margin=2.5cm]{geometry}
\usepackage{graphicx}
\usepackage{hyperref}
\usepackage{mathtools}
\usepackage{multirow}
\newcommand{\tensorTwo}[1]{\oalign{\mbox{$\bf{{#1}}$}\crcr\hidewidth$\scriptscriptstyle\sim$\hidewidth}}
\newcommand{\tensorFour}[1]{\tensorTwo{\tensorTwo{#1}\vphantom{#1}}}
\newcommand{\Stress}{\tensorTwo{\bm{\sigma}}}
\newcommand{\StressRate}{\tensorTwo{\dot{\bm{\sigma}}}}
\newcommand{\Strain}{\tensorTwo{\bm{\epsilon}}}
\newcommand{\StrainRate}{\tensorTwo{\dot{\bm{\epsilon}}}}
\renewcommand{\vec}{\bm}
\newcommand{\R}{\mathbb{R}}
\newcommand{\dx}{\;d\x}

\newcommand{\dz}{\;d\z}
\newcommand{\grad}{\nabla}
\newcommand{\tr}{^\top}
\newcommand{\inv}{^{-1}}
\newcommand{\defas}{\coloneqq}

\newcommand{\expct}[2]{\mathbb{E}_{#1}\left[#2\right]}

\DeclarePairedDelimiter{\norm}{\lVert}{\rVert}
\DeclarePairedDelimiter{\ceil}{\lceil}{\rceil}
\newcommand{\I}{\vec{I}}

\newcommand{\alpharesmin}{\alpha_{\text{res}}^\text{l}}
\newcommand{\alpharesmax}{\alpha_{\text{res}}^\text{h}}
\newcommand{\deltaalphares}{\Delta\alpha_{\text{res}}}
\newcommand{\lambdadotmax}{\dot{\lambda}^{\ast}}
\newcommand{\malpha}{m_{\alpha}}
\newcommand{\betamax}{\beta^{\ast}}
\newcommand{\lambdamax}{\lambda^{\ast}}
\newcommand{\mbeta}{m_{\beta}}
\newcommand{\strain}{\epsilon}
\newcommand{\stress}{\sigma}
\newcommand{\dataStrain}{\data_{\strain}}
\newcommand{\dataStress}{\data_{\stress}}
\newcommand{\dimStrain}{{n_\strain}}
\newcommand{\dimStress}{{n_\stress}}

\newcommand{\modelStrain}{\model_\strain}
\newcommand{\modelStress}{\model_\stress}
\newcommand{\covStrain}{{\cov_\strain}}
\newcommand{\covStress}{{\cov_\stress}}

\newcommand{\prior}{\rho_{\text{prior}}}
\newcommand{\posterior}{\rho_{\text{post}}}
\newcommand{\like}{\rho_{\text{like}}}
\newcommand{\data}{{\vec{d}}}
\newcommand{\dimData}{{n_\data}}
\newcommand{\dimPar}{n}
\renewcommand{\dim}{\dimPar}
\newcommand{\model}{\mathcal{G}}
\newcommand{\noise}{\vec{\eta}}
\newcommand{\cov}{\Gamma}
\newcommand{\misfit}{f_{\data}}
\newcommand{\x}{\vec{x}}

\newcommand{\C}{\vec{C}}
\newcommand{\W}{\vec{W}}

\newcommand{\MeigVals}{\Lambda}
\newcommand{\eigVal}{\lambda}
\newcommand{\eigVec}{\vec{w}}
\newcommand{\y}{{\vec{y}}}
\newcommand{\z}{{\vec{z}}}
\newcommand{\priory}{\rho_{\text{prior},\y}}
\newcommand{\priorzy}{\rho_{\text{prior},\z|\y}}
\newcommand{\prioryEst}{\hat{\rho}_{\text{prior},\y}}

\newcommand{\NxEss}{{N_{\x,\text{ESS}}}}
\newcommand{\NyEss}{{N_{\y,\text{ESS}}}}
\newcommand{\NylEss}{{N_{\y^{(l)},\text{ESS}}}}
\newcommand{\NzEss}{{N_{\z,\text{ESS}}}}
\newcommand{\NziEss}{{N_{\z,\text{ESS},i}}}
\newcommand{\Jmax}{J_\text{max}}

\makeatletter
\renewcommand{\email}[2][]{%
	\ifx\emails\@empty\relax\else{\g@addto@macro\emails{,\space}}\fi%
	\@ifnotempty{#1}{\g@addto@macro\emails{\textrm{(#1)}\space}}%
	\g@addto@macro\emails{#2}%
}
\makeatother
\renewcommand{\vec}[1]{\mathbf{{#1}}}
\graphicspath{{media/}}
\title{Efficient parameter estimation for a methane hydrate model with active subspaces}
\author[M. Teixeira Parente]{Mario Teixeira Parente$^\dagger{}^\P$}
\author[S. Mattis]{Steven Mattis$^\dagger$}
\author[S. Gupta]{Shubhangi Gupta$^*$}
\author[C. Deusner]{Christian Deusner$^*$}
\author[B. Wohlmuth]{Barbara Wohlmuth$^\dagger{}^\text{\textdaggerdbl}$}
\address{$^\dagger$ Chair for Numerical Mathematics, Technical University Munich, Germany}
\address{$^\text{\textdaggerdbl}$ Department of Mathematics, University of Bergen, Norway}
\address{$^*$ GEOMAR Helmholtz Centre for Ocean Research Kiel, Germany}
\address{$^\P$ \textnormal{Corresponding author}}
\email{$\lbrace \text{parente,mattis,wohlmuth}\rbrace$@ma.tum.de}
\email{$\lbrace \text{sgupta,cdeusner}\rbrace$@geomar.de}
\date{\today}
\begin{document}
\begin{abstract}
Methane gas hydrates have increasingly become a topic of interest because of their potential as a future energy resource.
There are significant economical and environmental risks associated with extraction from hydrate reservoirs, so a variety of multiphysics models have been developed to analyze prospective risks and benefits.
These models generally have a large number of empirical parameters which are not known a priori.
Traditional optimization-based parameter estimation frameworks may be ill-posed or computationally prohibitive.
Bayesian inference methods have increasingly been found effective for estimating parameters in complex geophysical systems.
These methods often are not viable in cases of computationally expensive models and high-dimensional parameter spaces.
Recently, methods have been developed to effectively reduce the dimension of Bayesian inverse problems by identifying low-dimensional structures that are most informed by data.
Active subspaces is one of the most generally applicable methods of performing this dimension reduction.
In this paper, Bayesian inference of the parameters of a state-of-the-art mathematical model for methane hydrates based on experimental data from a triaxial compression test with gas hydrate-bearing sand is performed in an efficient way by utilizing active subspaces.
Active subspaces are used to identify low-dimensional structure in the parameter space which is exploited by generating a cheap regression-based surrogate model and implementing a modified Markov chain Monte Carlo algorithm.
Posterior densities having means that match the experimental data are approximated in a computationally efficient way.
\end{abstract}
\maketitle
\section{Introduction}
\label{sec:intro}

Methane gas hydrates are crystalline solids formed when water molecules enclatharate methane molecules \cite{Sloan1998}.
Gas hydrates are stable at low temperatures and high pressures and occur naturally in permafrost regions and marine off-shores \cite{DaweThomas2007}.
If warmed or depressurized, gas hydrates destabilize and dissociate into water and gas.
It is estimated that the energy content of methane occurring as hydrates exceeds the \textit{combined} energy content of all other conventional fossil fuels \cite{Pinero2013}.
Natural gas hydrates are, therefore, deemed a promising future energy resource.
Several methods have been proposed for gas extraction from hydrate reservoirs, such as thermal stimulation, depressurization, and chemical activation \cite{Moridis2009,Moridis2011}. 
Application of these methods at large scales is, however, very challenging due to the inherent geo\-tech\-nic\-al risks associated with gas hydrate destabilization,
such as, rapid consolidation, seafloor subsidence, well collapse, uncontrolled sand migration, and local and regional slope instability \cite{Sultan2004a,Sultan2004b}.
In order to quantify these risks under various production scenarios and to make realistic assessments regarding the viability of these production methods,
a number of multiphysics models (e.g., \cite{SGupta2015,Hyodo2014,Kimoto2010,KlarSogaNG2010,KlarUchidaSogaYamamoto2013,Rutqvist2011}) have been developed in the recent years.
It is known that the gas hydrate-bearing sediments (GHBS) are very complex geo\-mat\-er\-i\-als which show a wide range of geomechanical behaviors depending on their distribution, saturation, morphology, formation, and consolidation history, etc.
The predictive capability of these models, therefore, depends heavily on the accuracy of the input constitutive model and parameters. 
A number of constitutive models have been proposed so far to describe the geomechanical behavior of GHBS \cite{Pinkert2017,PinkertGrozicPriest2015,UchidaSogaYamamoto2012MHCS,Sanchez2017}.
One common feature of these models is their large number of empirical parameters, often exceeding ten.
The models themselves are highly complex, and traditional techniques of estimating the model parameters not only require large experimentation effort, but also very large (often, even prohibitive) computational efforts from solving multi-dimensional nonlinear optimization problems which may be ill-posed.

In recent years, Bayesian methods incorporating model and data uncertainties have been successful for constructing well-posed, solvable parameter estimation problems.
Bayesian inference and Generalized Likelihood Uncertainty Estimation (GLUE) have proven to be among the most popular approaches for quantifying uncertainties in problems in porous media; e.g. \cite{UQ-GLUE_Beven_Freer,Bayesian_GLUE_Freer,Bayesian_Leube,Bayesian_Nowak,Bayesian_Troldborg,GLUE_Vrugt}. 
In these approaches, the misfit between experimental data and evaluations of the mathematical model is used to define a statistical map called the likelihood function.
When a prior distribution is defined on the parameters (incorporating prior knowledge of the physics and the model), a posterior distribution is defined in terms of the likelihood and prior distribution using Bayes' Theorem.
A desired property of a posterior distribution is that when it is propagated through the forward map, it matches well with the measured data.
In this paper, we analyze how well the forward mapping of the mean of the posterior matches with the experimental data; however, there are more recent notions of consistency of solutions of Bayesian inverse problems\cite{butler2018combining}, which could be utilized in future work.
Common objectives are to produce a set of samples following the distribution of the posterior or to determine the parameters of maximum likelihood.
These objectives are often achieved using a Markov chain Monte Carlo (MCMC) method.

MCMC methods can be computationally prohibitive if the forward model is computationally expensive because it often must be evaluated a large number of times.
Additional computational complexity occurs if the space of uncertain parameters is high-dimensional because Markov chains must explore the parameter space to find regions of high probability.
There recently has been much effort in reducing the computational expense of MCMC in such settings by exploiting the structure of the stochastic inverse problem with adaptive sampling methods \cite{haario2006dram,vrugt2009accelerating} and methods that effectively reduce the dimension of the parameter space \cite{bui2012extreme,bui2014solving,constantine2016accelerating,cui2016dimension,martin2012stochastic}.
The likelihood-informed subspace method \cite{cui2016dimension} and active subspace method \cite{constantine2016accelerating} identify data-informed subspaces which can be utilized to accelerate MCMC.
Likelihood-informed subspaces have restrictive conditions on the prior that are not met by many problems.
MCMC with active subspaces has more general conditions and is applicable to any problem where an active subspace exists.

Little effort has been put forth in the accurate estimation of parameters for state-of-the-art methane hydrate models because the high computational cost of model evaluations and the high-di\-men\-sio\-nal\-i\-ty of the parameter space make traditional methods computationally prohibitive.
However, relatively little is known about the relationships between the parameters in these models, and there are likely lower-dimensional structures in the parameter space that have not yet been identified which may be used to effectively reduce the dimension of the inverse problem.
We hypothesize that the method of active subspaces can be used to identify such lower-dimensional structures within the space of model parameters and that MCMC with active subspace can be used to efficiently perform Bayesian inference on the parameters which would otherwise be extremely computationally expensive.

The paper is organized as follows.
In Section \ref{sec:experiment}, details of the experimental study from which the data is generated are presented.
A mathematical model for deformation of a porous medium with methane hydrates is developed in Section \ref{sec:model}.
Section \ref{sec:inversion} presents a framework for Bayesian inversion with active subspaces for the methane hydrate problem.
Results from the inversion framework applied to the model using experimental data are shown in Section \ref{sec:results}, and conclusions are discussed in Section \ref{sec:conclusions}.

\section{Experimental study}
\label{sec:experiment}

Experimental data were obtained in a controlled triaxial compression test with gas hydrate-bearing sand (GHBS). 
GHBS was formed under controlled isotropic effective stress using the \textit{excess-gas-method} \cite{Choi2014,Priest2009}. 
In the excess-gas-method, gas hydrates are formed in partially water saturated porous or granular media by supplying gas within gas hydrate stability boundaries, i.e. at high pressure and low temperature.   
The use of the excess-gas-method enables the formation of homogeneously distributed gas hydrates in the porous matrix and adjustment of well-defined gas hydrate saturations ($S_h$) as a consequence of the limited availability of water.
Further, due to initial phase distributions and wetting behavior, gas hydrates are preferentially formed on grain surfaces and in pore throats. 
This microscale phase distribution is recognized to result in mechanical strengthening of the bulk sediment \cite{Hyodo2005,Masui2005}.    
After completion of methane hydrate formation, drained triaxial compression tests were performed at controlled axial strain rates under quasi-static loading and constant confining effective stress. 

\subsection{Experimental setup and measurements}\label{subsec:setup_and_components} 

\begin{figure}
	\centering
	\includegraphics[scale=0.5]{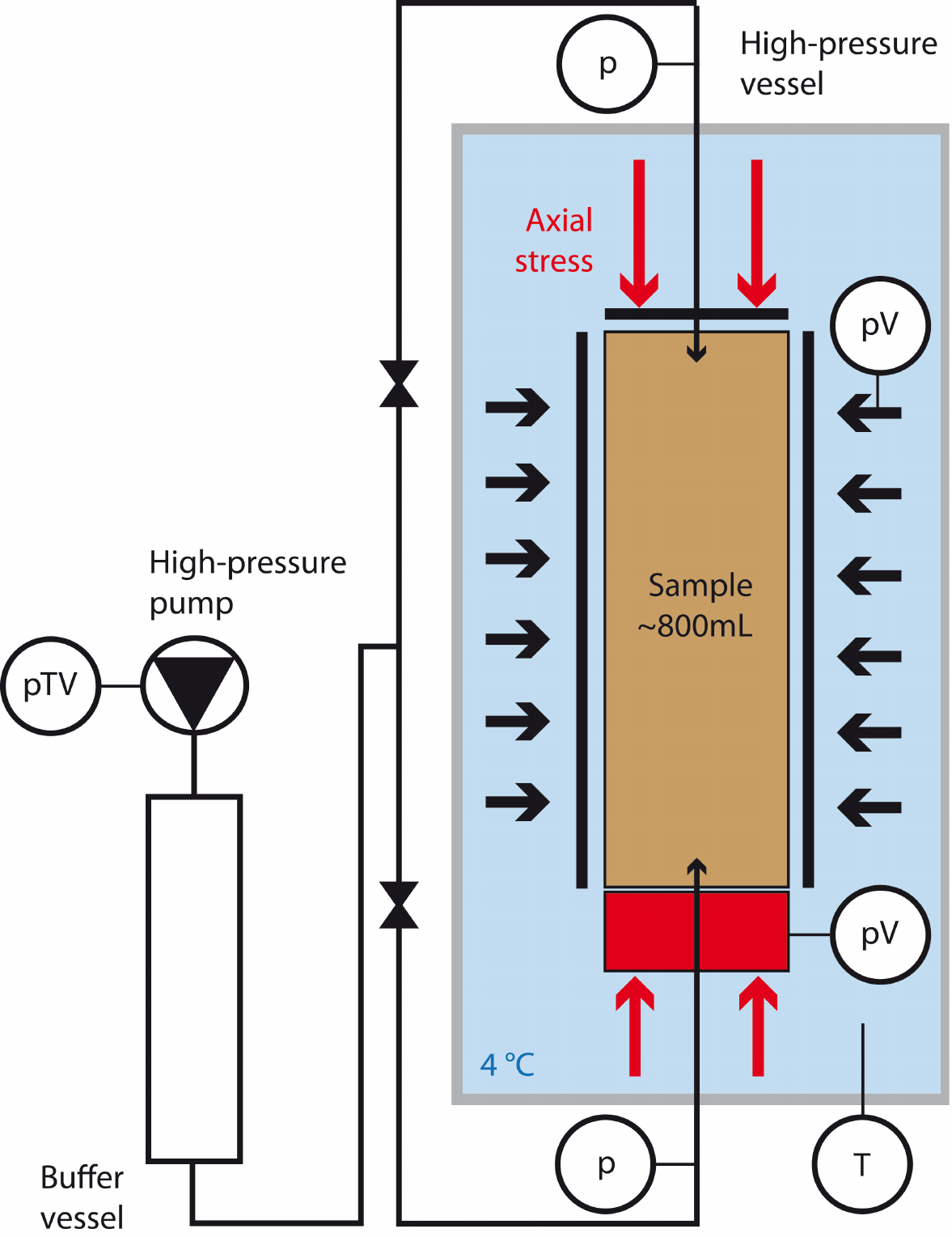}
	\caption{Experimental scheme and measurements.}
	\label{fig:experimental_setup}
\end{figure}

Experiments were carried out in the custom-made high pressure apparatus NESSI (Natural Environment Simulator for Sub-seafloor Interactions) \cite{Deusner2012} (see Fig.~\ref{fig:experimental_setup}), 
which is equipped with a triaxial cell mounted in a $40$ l stainless steel vessel (APS GmbH Wille Geotechnik, Rosdorf, Germany). 
The sample sleeve is made from FKM. Other wetted parts of the setup are made of stainless steel. 
Axial and confining stresses and sample volume changes were monitored throughout the overall experimental period using high-precision hydraulic pumps. 
Pore pressure was measured in the influent and the effluent fluid streams close to the sample top and bottom. 
Pressure control was achieved using automated high-pressure piston pumps (Teledyne ISCO, Lincoln, USA). 
The experiment was carried out under constant temperature conditions, temperature control was achieved with a thermostat system (T1200, Lauda, Lauda-K\"onigshofen, Germany). 

\label{subsec:measurements}
Experimental control and process monitoring was carried out using high-precision piston pumps which individually control pressure and volumes of hydraulic (axial and confining) and pore fluids. 
During triaxial compression, pressure and fluid volumes were measured and recorded at $1$~s intervals. 
The accuracy of the individual pressure measurements is $\pm  0.5\%$ at constant temperature. 
Random errors resulting from temperature changes or leakage of fluids can be neglected due to the short duration of the compression tests and large thermal buffer capacity of the high-pressure systems. 
The accuracy of volume and strain measurements is related to pressure measurements since system volume changes are calibrated depending on the system pressure. 
Thus, erroneous pressure measurements can result in an overall error of volume measurement of $4$ ml, which converts to $0.4\%$ of volumetric strain.          

The sediment sample was prepared from quartz sand (initial sample porosity: $0.35$, grain size: $0.1 - 0.6$ mm, G20TEAS, Schlingmeier, Schw\"ulper, Germany),
which was mixed with de-ionized water to achieve a final water saturation of $0.2$ relative to the initial sample porosity. 
The partially water-saturated and thoroughly homogenized sediment was filled into the triaxial sample cell to obtain final sample dimensions of $160$ mm in height and $80$ mm in diameter. 
The sample geometry was assured using a sample forming device. 
The sample was cooled to $4$ \degree C after the triaxial cell was mounted inside the pressure vessel. 

\subsection{Experimental procedure}

Prior to the gas hydrate formation, the partially water-saturated sediment sample was isotropically consolidated to $1$ MPa effective stress under drained conditions. 
The sample was flushed with $CH_4$ gas and, subsequently, pressurized with $CH_4$ gas to obtain a pore pressure of approximately $10$ MPa. 
During pressurization with $CH_4$ gas, and throughout the overall gas hydrate formation period, isotropic effective stress was controlled to remain constant at $1$ MPa using an automated control algorithm. 
The formation process was continuously monitored by logging the $CH_4$ gas pressures. 
Mass balances and volume saturations were calculated based on $CH_4$ gas pressure to confirm that available pore water was fully converted into gas hydrates.  

After completion of gas hydrate formation, the triaxial compression test was conducted at a controlled axial strain rate of $0.1$ mm/min. 
During axial loading and compression, the confining effective stress was controlled to remain constant by adjusting the confining hydraulic fluid volume in the pressure vessel.
Accumulated volumetric strain was calculated based on changes of axial and confining volumes, which are monitored by the hydraulic pumps.    

\section{Methane hydrate model}
\label{sec:model}

From a geomechanics point of view, the gas hydrate-bearing sands (GHBS) are cohesive-frictional granular materials.
The most important features of the mechanical behavior of GHBS include plastic deformations and the onset of critical state (i.e. i\-so\-cho\-ric deformations),
pressure-dependence, and shear-volumetric strain coupling (or \textit{dilatancy}).
The strength and the stiffness of GHBS are strongly influenced by gas hydrate saturation and hydrate morphology (i.e. pore-filling, load-bearing, cementing, etc.), 
as well as the hydrate formation method \cite{WaiteSantamarina2009,Yun2007}.
In general, the higher the gas hydrate saturation, the higher is the bulk compressive strength of the GHBS.
The gas hydrate saturation also enhances the cohesive strength, frictional resistance, and dilatancy of GHBS. 
We model the geomechanical behavior of GHBS within an incremental elasto-plasticity framework, 
and for simplicity we assume infinitesimal strains.
This section presents the main elements of our elasto-plastic material model for GHBS, including the yield function, plastic flow directions, and hardening and softening evolution laws.

\paragraph{Notation}
For any second order tensors $\tensorTwo{a}$ and $\tensorTwo{b}$, the inner product is given by $\tensorTwo{a} : \tensorTwo{b} = a_{ij} b_{ij} $,
and the dyadic product is given by $\left( \mathbf{a} \otimes \mathbf{b}\right) = a_{ij} b_{kl}$.
The tensor product between a fourth order tensor $\tensorFour{A}$ and a second order tensor $\tensorTwo{b}$ is given by $ \tensorFour{A} : \tensorTwo{b} = A_{ijkl} b_{kl}$.
The Euclidean norm of $\tensorTwo{a}$ is given by $\lVert \tensorTwo{a} \rVert = \left(\tensorTwo{a}:\tensorTwo{a}\right)^{1/2}$.
The second order unit tensor is given by $\tensorTwo{I}=\delta_{ij}$, where $\delta$ denotes the Kronecker delta function.
The fourth order unit tensor is given by $\tensorFour{I} = \frac{1}{2}\left(\delta_{ik}\delta_{jl}+\delta_{il}\delta_{jk}\right)$.
The trace of $\tensorTwo{a}$ is given by $Tr\tensorTwo{a} = \tensorTwo{I}:\tensorTwo{a}$.
Any second order tensor $\tensorTwo{a}$ can be decomposed into a dilational (or volumetric) part, dil$\tensorTwo{a} $, and a deviatoric part, dev$\tensorTwo{a} $.
In $3$D, 
dil$\tensorTwo{a} = \frac{1}{3} Tr \tensorTwo{a} $ and
dev$\tensorTwo{a} = \tensorTwo{a} - \frac{1}{3} Tr \tensorTwo{a} $ .

\subsection{Preliminaries}
\label{subsec:Preliminaries}

Let $\Stress$ be the Cauchy stress tensor and $\Strain = \frac{1}{2}\left( \nabla \mathbf{u} + \nabla\tr \mathbf{u} \right)$ the infinitesimal strain tensor.
The vector $\mathbf{u}$ denotes the displacement field.
Both $\Stress$ and $\Strain$ are symmetric second order tensors.
The total infinitesimal strain $\Strain$ is decomposed \textit{additively} into the elastic strain $\Strain^e$ and the plastic strain $\Strain^p$, i.e.
$
\Strain = \Strain^e + \Strain^p.
$

In classical plasticity \cite{jirasekBazant2002}, the state of stress depends on the loading-unloading history and is calculated incrementally. 
The stress and the strain rate tensors (i.e. $\StressRate$, $\StrainRate$, $\StrainRate^e$, and $\StrainRate^p$) are approximated using an implicit Euler finite difference method, i.e. for any time interval $\left[ t_i , t_{i+1} \right]$, the stress or strain rates are approximated as 
\begin{equation}
\tensorTwo{\dot{\left[\cdot\right]}} = \dfrac{\tensorTwo{\left[\cdot\right]} - \tensorTwo{\bar{\left[\cdot\right]}}}{t_{i+1}-t_i},
\end{equation}
where $\tensorTwo{\bar{\left[\cdot\right]}}$ denotes the state at time $t_i$, and $\tensorTwo{\left[\cdot\right]}$ denotes the state at time $t_{i+1}$.

We define the plasticity relationships in terms of the following stress invariants:
\begin{align}
 \label{eqn:StressInvariants}
 p = \frac{1}{3} Tr \Stress \quad \text{and} \quad q = \sqrt{\frac{3}{2}} \lVert \text{dev } \Stress \rVert,
\end{align}
where $p$ denotes hydrostatic or mean stress, and $q$ denotes shear stress.
Corresponding invariants of strain rate are:
\begin{align}
 \label{eqn:StrainInvariants}
 \dot\epsilon_v = Tr \StrainRate \quad \text{and} \quad \dot\epsilon_s = \sqrt{\frac{2}{3}} \lVert \text{dev } \StrainRate \rVert,
\end{align}
where $\epsilon_v$ denotes the volumetric strain rate, and $\dot\epsilon_s$ denotes the shear strain rate.
The invariants of elastic and plastic strain rates are defined similarly.

\subsection{Elasticity}
\label{subsec:Elasticity}

In the elastic range, we assume linear isotropic material behavior, i.e. the stress $\Stress$ is related to the elastic strain $\Strain_e$ through a linear Hooke's law,
\begin{align}
 \label{eqn:ElasticHooksLaw}
 \Stress = \tensorFour{C}^e : \Strain^e,
\end{align}
where $\tensorFour{C}^e$ is the elastic stiffness of the material which is a symmetric positive definite fourth order tensor,
\begin{align}
 \label{eqn:ElasticStiffnessTensor}
 \tensorFour{C}^e = L_1 \tensorTwo{I} \otimes \tensorTwo{I} + 2 L_2 \tensorFour{I}.
\end{align}
$L_1$ and $L_2$ are the Lam\'e coefficients,
\begin{align}
 \label{eqn:LameCoefficients}
 L_1 = \frac{\nu E\left(S_h\right)}{\left(1+\nu\right) \left(1-2\nu\right)} \quad \text{and} \quad L_2 = \frac{E\left(S_h\right)}{2\left(1+\nu\right)},
\end{align}
where $E$ is the elastic Young's modulus, and $\nu$ is the Poisson's ratio of the GHBS, respectively.
$S_h$ denotes the gas hydrate saturation.
It is observed that the Young's modulus of GHBS increases with increasing $S_h$,
while the Poisson's ratio does not vary much over a wide range of $S_h$ and can be assumed constant \cite{Lee2010,Miyazaki2011}.
In \cite{SantamarinaRuppel2010}, the authors have proposed an empirical relationship for $E$ of the form
\begin{align}
 \label{eqn:YoungsModulous}
 E = E_s\left(\sigma_c\right) + S_h^{m} E_h,
\end{align}
where $E_s$ and $E_h$ denote the Young's modulus of the sand and gas hydrates, respectively, and $\sigma_c$ is the confining stress.
The exponent $m$ varies over a wide range. 
In their experiment and modeling study in \cite{Gupta2017}, the authors found that the effect of $S_h$ on $E$ was linear ($m=1$) during hydrate formation, 
while during hydrate dissociation the effect of $S_h$ on $E$ was stronger ($m=3$).

\subsection{Yield function}
\label{subsec:YieldFunction}

There exists a yield surface $F$ in the stress space that encompasses the elastic region. 
The stress states lying \textit{inside} the yield surface produce elastic deformations, while the stress states lying \textit{on} the surface produce plastic deformations.
The stress states outside the yield surface are inadmissible.

\begin{figure}
	\centering
	\includegraphics[scale=0.5]{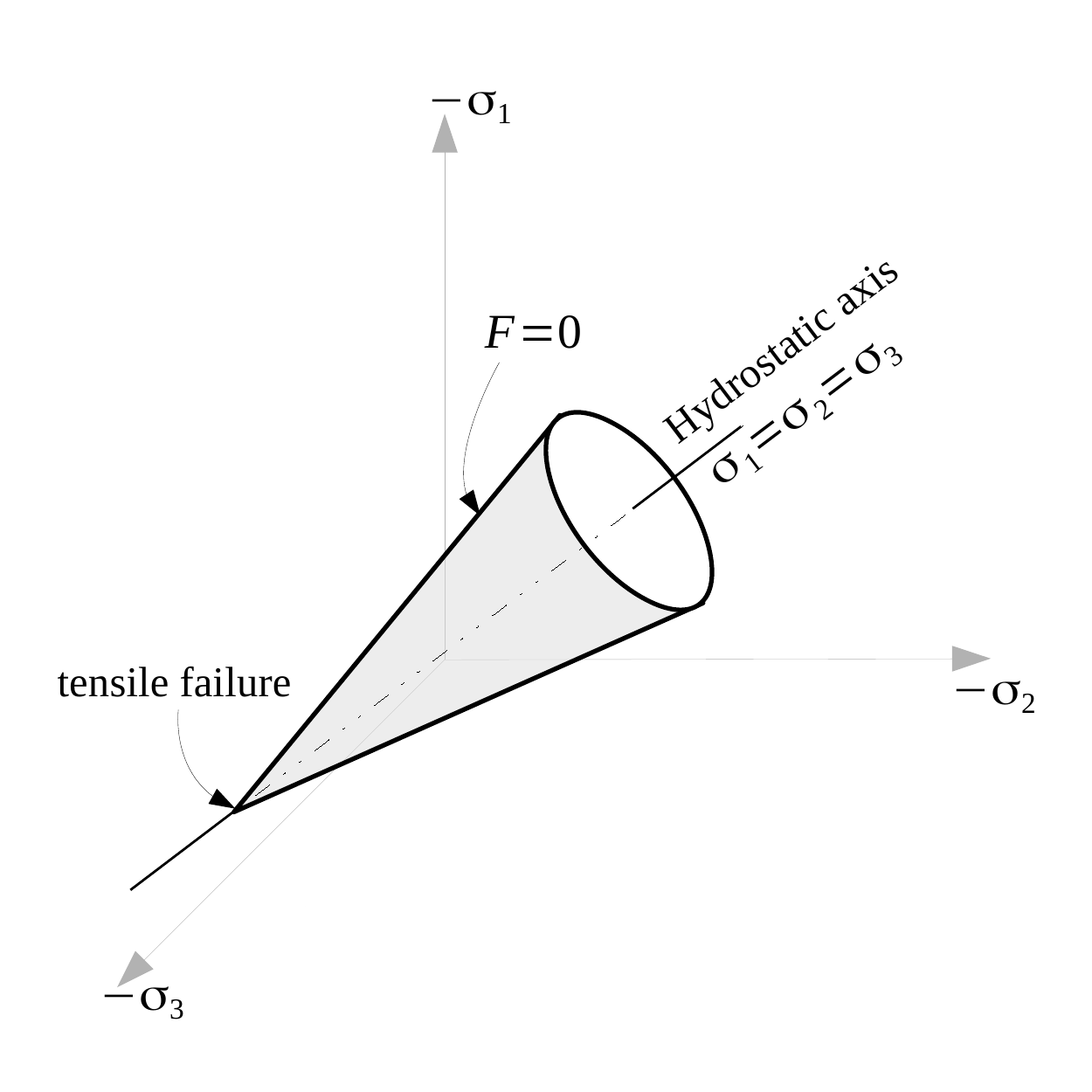}
	\caption{Drucker-Prager yield surface in principal stress space.}
	\label{fig:DPYieldFunctionInPrincipalStressSpace}
\end{figure}

We consider a Drucker-Prager yield criterion where the yield function is given as
\begin{align}
 \label{eqn:YieldFunction}
 F\left( \Stress,\bm{\chi} \right) \defas q + \alpha\left( \bm{\chi} \right) p - c\left( \bm{\chi} \right) = 0.
\end{align}
Function~$F$ describes a conical surface in the principal stress space (see Fig.~\ref{fig:DPYieldFunctionInPrincipalStressSpace}).
The parameter $\alpha$ indicates the mobilized frictional resistance at any given stress state.
The parameter $c$ indicates the cohesive strength of the granular material.
$\bm{\chi}$ denotes the vector of internal plastic variables which affect the hardening-softening behavior of GHBS 
due to changes in internal structure or grain contacts, packing density of the sand grains, hydrate saturation, hydrate pore habit, etc.

\subsection{Plastic strains}
\label{subsec:PlasticStrains}

\begin{figure}
	\centering
	\includegraphics[scale=0.5]{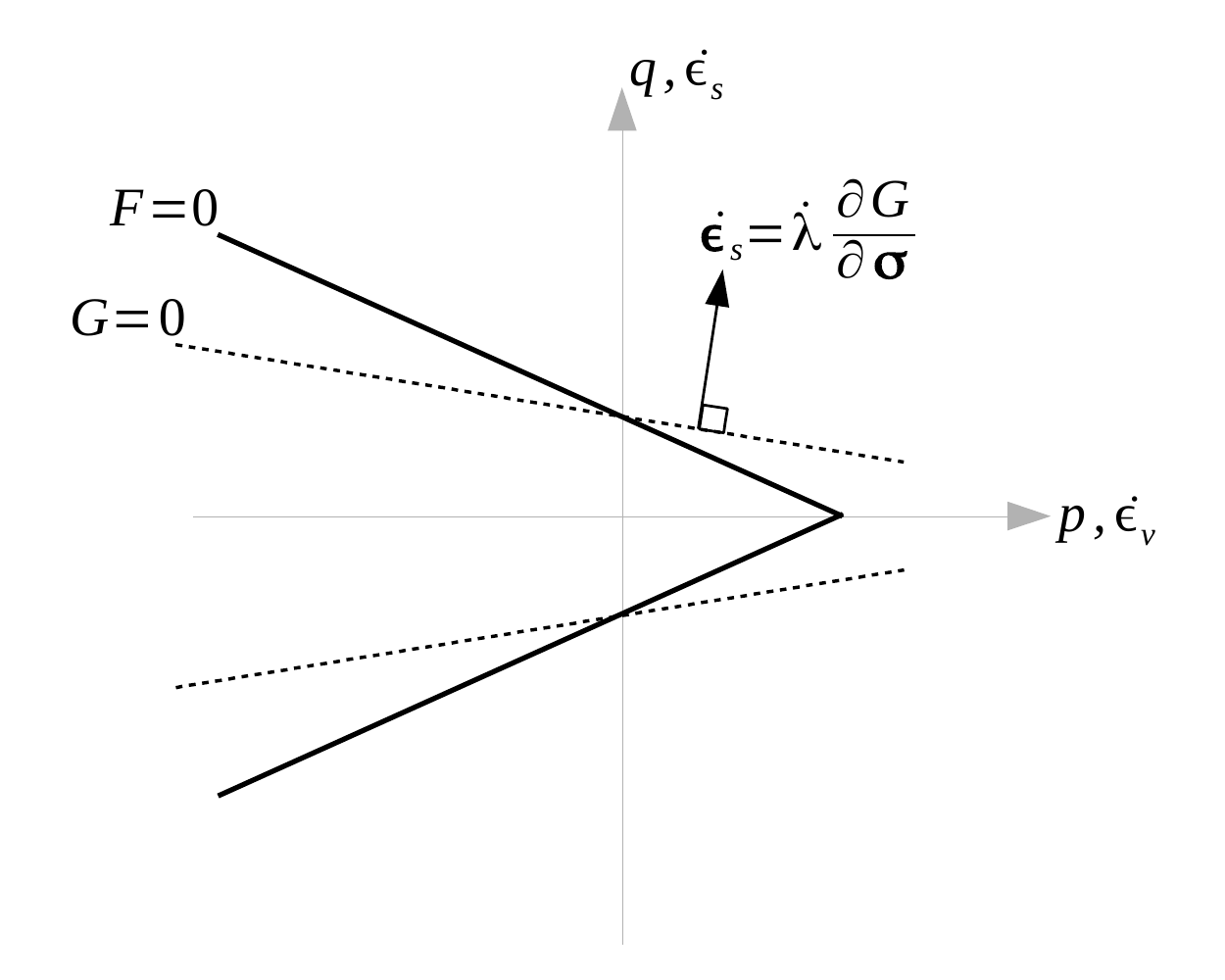}
	\caption{Potential surface and plastic strains in $p-q$ space.}
	\label{fig:PlasticStrains}
\end{figure}

Similar to the yield surface  $F$, there exists a plastic potential surface $G$ in the stress space such that the plastic flow occurs in a direction normal to this surface (see Fig.~\ref{fig:PlasticStrains}).
The incremental plastic strains (i.e. the plastic strain rate) can be derived from the plastic potential $G$ as
\begin{align}
 \label{eqn:PlasticStrains}
 \StrainRate^p = \dot{\lambda} \frac{\partial G}{\partial \Stress},
\end{align}
where $\partial G/\partial \Stress$ describes the normal to the surface $G$, and $\dot{\lambda}$ is a proportionality constant indicative of the magnitude of the plastic strain increment.
It can further be shown that the invariants of the plastic strain rate can be written as
\begin{align}
 \label{eqn:PlasticStrainInvariants}
 \dot{\epsilon}^p_v = \dot{\lambda} \frac{\partial G}{\partial p} 
 \quad \text{and} \quad
 \dot{\epsilon}^p_s = \dot{\lambda} \frac{\partial G}{\partial q}.
\end{align}

We consider a non-associative flow rule, i.e. $G \neq F$, 
\begin{align}
 \label{eqn:PlasticPotentialSurface}
 G\left( \Stress,\bm{\chi} \right) \defas q + \beta\left( \bm{\chi} \right) p = 0,
\end{align}
where $\beta < \alpha$.
The parameter $\beta$ denotes the \textit{dilatancy} of the material. 
Dilatancy is a characteristic property of frictional granular materials. 
It contributes to the strength of the material and effectively couples the deviatoric and volumetric components of plastic deformation.
From Eqns.~\eqref{eqn:PlasticStrainInvariants} and \eqref{eqn:PlasticPotentialSurface}, dilatancy can be written as $\beta = \dot{\epsilon}^p_v/\dot{\epsilon}^p_s$.
Depending on the relative packing density of the grains, it can allow for macroscopic contraction or dilation of the material under external loads.

\subsection{Loading-unloading conditions}
\label{subsec:LoadingUnloadingConditions}

Along any process of loading-unloading, if $F < 0$, the stress state is elastic and $\dot{\lambda} = 0$,
while, if $F=0$, the stress state is plastic and $\dot{\lambda} > 0$.
These nonlinear inequality constraints can be reformulated as the following Karush-Kuhn-Tucker \cite{KuhnTucker1951} optimality conditions:
\begin{align}
\label{eqn:KuhnTuckerConditions}
 F\left( \Stress,\bm{\chi} \right) \leq 0,\quad \dot{\lambda} \geq 0, \quad \dot{\lambda} F = 0.
\end{align}
To confine the stress trajectory to the yield surface during plastic loading, an additional plastic consistency condition is considered \cite{Owen2011}:
\begin{align}
 \label{eqn:PlasticConsistencyCondition}
 \dot{\lambda} \dot{F} = \dot{\lambda} \left( \frac{\partial F}{\partial \Stress} : \StressRate + \frac{\partial F}{\partial \bm{\chi}} : \dot{\bm{\chi}} \right) = 0.
\end{align}

\subsection{Evolution laws}
\label{subsec:EvolutionLaws}

Experiments have shown that an increase in gas hydrate saturation tends to increase the initial frictional resistance, apparent cohesive resistance, peak strength, and peak dilatancy of GHBS 
\cite{Hyodo2005,Lee2007,Masui2005}.
It is also observed that GHBS show a distinct strain hardening-softening behavior \cite{Miyazaki2011,Miyazaki2010}.
In \cite{UchidaSogaYamamoto2012MHCS}, this behavior is explained in detail, and strain dependent evolution laws are presented to capture the observed strain hardening-softening.
Other strain dependent evolution laws have also been presented in \cite{KlarSogaNG2010,PinkertGrozic2014,PinkertGrozicPriest2015}, among others.
In our experiments, in addition to strain hardening-softening, we additionally observe a distinct \textit{secondary} hardening phase in the stress-strain response of the GHBS samples.

Assuming that the frictional resistance of any geomaterial can be expressed as a sum of dilatancy and some \textit{residual} frictional resistance \cite{Wood1990}, i.e.
$$ \alpha = \beta + \alpha_{\text{res}}$$
where, at critical state, $\beta = 0$ and $\alpha = \alpha_{\text{res}}$,
we hypothesize that the primary hardening occurs due to the dilatancy of the sample, 
while the secondary hardening occurs due to an increase in residual frictional resistance under plastic loading.
This is likely the result of internal damage and hydrate redistribution in the pore spaces, causing higher particle density, increased interlocking of grains, and higher friction at the grain contacts.
This effect becomes dominant when the material has spent all its dilatancy and has achieved a critical state.
We ignore the contribution of the cohesive strength on the observed strain dependent hardening-softening-hardening behavior.

We consider the plastic internal variables $\bm{\chi} = \left( \dot{\epsilon}^p_s ,\epsilon^p_s, S_h \right)\tr$.
From Eqns. \eqref{eqn:PlasticStrainInvariants} and \eqref{eqn:PlasticPotentialSurface}, we get $\dot{\epsilon}^p_s = \dot{\lambda}$ and $\epsilon^p_s = \int^{t_{i+1}}_{t_i}\dot{\lambda} \;dt = \lambda$.
To capture the macroscopic stress-dilatancy behavior of GHBS observed during our triaxial compression experiments,
we describe smooth empirical evolution laws for the properties $\alpha$, $\beta$, and $c$ in Eqns. \eqref{eqn:YieldFunction} and \eqref{eqn:PlasticPotentialSurface} as
\begin{align}
 \label{eqn:CohesionEvolutionLaw}
 &c = c\left(S_h\right),
 \\
 \label{eqn:DilatancyEvolutionLaw}
 &\beta = \beta^{\ast}\left(S_h\right)
	 \cdot \bar{\lambda}
         \cdot \exp\left( 1 - \bar{\lambda}^{m_{\beta}} \right) ,
 \\
 \label{eqn:FrictionalResistanceEvolutionLaw}
 &\alpha_{\text{res}} = \alpharesmin\left(S_h\right) + \deltaalphares\left(S_h\right)
				 \cdot \left( 1+1/\dot{\bar{\lambda}}\right)^{-1}
				 \cdot \bar{\lambda}^{m_{\alpha}},
\end{align}
where $\bar{\lambda} = \lambda/\lambda^{\ast}\left(S_h\right)$ and $\dot{\bar{\lambda}} = \dot{\lambda}/\dot{\lambda}^{\ast}\left(S_h\right)$.
The functional dependence of the parameters $c$, $\beta^{\ast}$, $\lambda^{\ast}$, $\alpharesmin$, $\deltaalphares$, and $\dot{\lambda}^{\ast}$ on $S_h$ 
can be derived through empirical correlation by repeating these experiments over a range of hydrate saturations.
In this work, we consider only a single GHBS sample with a constant hydrate saturation.
So, the exact functional dependence of the plasticity parameters on $S_h$ is not of direct relevance for the purpose of presenting our use of the active subspace strategy and will not be discussed further.

Eqns. \eqref{eqn:DilatancyEvolutionLaw} and \eqref{eqn:FrictionalResistanceEvolutionLaw} are extensions of the evolution functions proposed in \cite{Andrade2012}. 
The parameter $\beta^{\ast}$ denotes the peak dilatancy, and $\lambda^{\ast}$ denotes the corresponding accumulated plastic shear strain.
The parameter $\alpharesmin$ denotes the minimum frictional resistance of the intact material before loading.
It is interesting to note that in Eqn. \eqref{eqn:FrictionalResistanceEvolutionLaw}, when $\dot{\epsilon}_s \rightarrow 0$, we get $\alpha_{\text{res}} = \alpharesmin$.
Physically, this implies that under quasi-static loading conditions, the material does not undergo microscopic damage, and the residual frictional resistance of the material remains constant.
At higher loading rates ($\dot{\epsilon}_s>0$), however, the effects of microscopic damage, sand and hydrate grain rearrangement, friction at grain contacts, etc. become progressively larger,
resulting in an overall increase in the macroscopic residual frictional resistance.

\subsection{Finite element implementation}
\label{subsec:FEImplementation}

We solve the \textit{global} nonlinear equilibrium equation using a Galerkin finite element formulation defined on $Q_1$ elements.
The nonlinearities are resolved iteratively using a full Newton-Raphson method with a continuum tangent matrix \cite{zienkiewiczTaylor2013}.
Within each global Newton iteration step, a \textit{local} problem is solved at each Gauss point to determine the new stress state.
The local problem involves the integration of the material model (Section~\ref{subsec:Elasticity}-\ref{subsec:EvolutionLaws}) over the load increment of the current global step.
We use an implicit return mapping algorithm \cite{HuangGriffiths2009,SimoHughes2006} to solve the local problem.
The implicit algorithm uses the \textit{final} point in the stress space to evaluate the relevant derivatives and variables.
Since this point is not known in advance, a Newton-Raphson method is used to advance the solution iteratively toward the final solution.
In a more generalized solution method,
the nonlinear equilibrium equation as well as the inequality constraints can be treated within a single Newton iteration,
which can be implemented as a primal-dual active set strategy (e.g. \cite{WohlmuthHager2009,WohlmuthHager2010}).
We have implemented our numerical scheme in C++ based on the DUNE PDELab framework \cite{Bastian2010,DUNEPdelab2012}.

\subsection{Numerical simulation of the triaxial compression experiments}
\label{subsec:NumericalSimulation}

\begin{figure}
	\centering
	\includegraphics[scale=0.5]{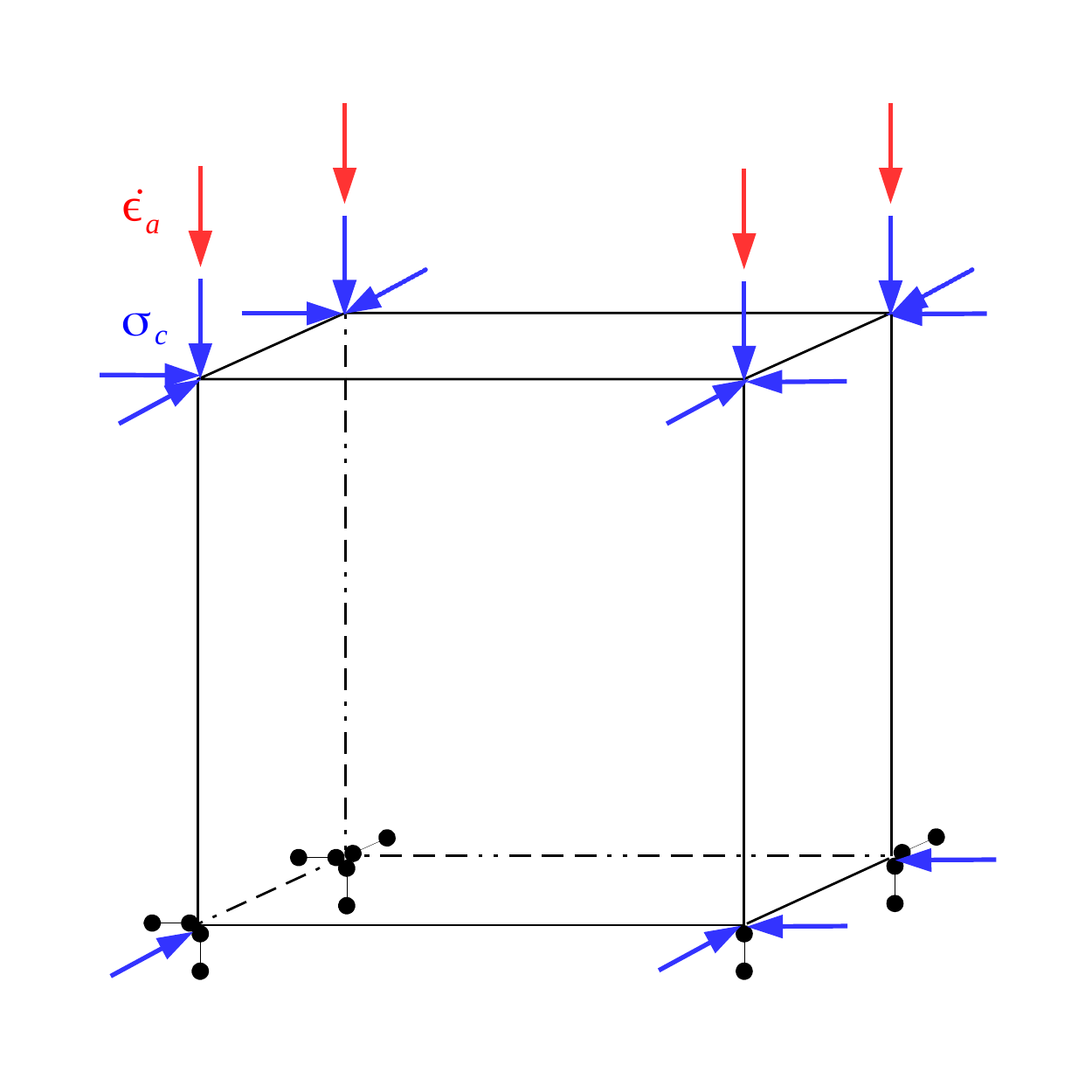}
	\caption{One element triaxial setup.}
	\label{fig:OneElementTriaxSchematic}
\end{figure}

We consider a \textit{one element} triaxial setup with unit dimensions, as shown in Fig.~\ref{fig:OneElementTriaxSchematic}.
Load is applied in two stages. 
In the first stage, an isotropic load equal to the confining stress of $\sigma_c=1$ MPa is applied. 
This corresponds to the initial stress state of the GHBS sample in the triaxial compression experiment. 
In the second stage, the strain-controlled triaxial compression of the GHBS sample is simulated 
by specifying an axial strain rate of $\dot{\epsilon}_a=-1.04167 \times 10^{-5}$ $\text{s}^{-1}$, 
which corresponds to a vertical displacement rate of $0.1$ mm/minute in the $-z$ direction.
The strain-controlled load is applied incrementally in $1350$ steps with a step-size of $10$ s.
\section{Bayesian inversion with active subspaces}
\label{sec:inversion}
This section explains an approach to efficiently infer the parameters
\begin{equation}
\x = (c, \alpharesmin, \deltaalphares, \lambdadotmax, \malpha, \betamax, \lambdamax, \mbeta)\tr
\end{equation}
of the model described in Section~\ref{sec:model}. Note that the parameter $\x^{(3)}=\deltaalphares$ is the difference $\deltaalphares=\alpharesmax-\alpharesmin$.
\subsection{Bayesian inversion}
The statistical inference of parameters is done here with a Bayesian approach to inverse problems \cite{stuart2010inverse}.
In this approach, the goal is to interrogate a probability measure on the space of parameters, incorporating prior knowledge and data.
Prior knowledge about the parameters from physics or engineering expertise are considered in a \textit{prior} probability density function $\prior$.
This knowledge is updated using the data to get a \textit{posterior} density function $\posterior(\cdot|\data)$.
The posterior is a probability density function on the parameter space that is conditioned on the data, $\data$.
The forward model is represented by a \textit{quantity of interest map} $\model:\R^\dim\to\R^\dimData$, which is related to the data by an additive noise model
\begin{equation}
	\data = \model(\x)  + \noise,
\end{equation}
where $\noise\sim\mathcal{N}(\vec{0}, \cov)$ is zero-centered Gaussian noise with covariance  $\cov$.
The posterior density can be investigated by Bayes' formula
\begin{equation}
	\label{eq:bayes}
	\posterior(\x|\data) = \frac{\like(\data|\x)\prior(\x)}{\int{\like(\data|\x')\prior(\x'){\dx'}}},
\end{equation}
where $\like$ is the \textit{likelihood function} given by $\like(\data|\x) \propto \exp(-\misfit(\x))$ for the \textit{data misfit function}
\begin{equation}
	\label{eq:misfit}
	\misfit(\x) \defas \frac{\norm{\data-\model(\x)}^2_\cov}{2} \defas \frac{\norm{\cov^{-1/2}\left(\data-\model(\x)\right)}_2^2}{2}.
\end{equation}

MCMC is a popular technique to sample from the posterior density because of its general applicability \cite{brooks2011handbook,kaipio2006statistical}.
The \textit{Metropolis-Hastings algorithm} is a classical algorithm for MCMC, which constructs a Markov chain that has a stationary distribution equal to the posterior.
One can neglect the normalizing constant in Eqn.~\eqref{eq:bayes} which often is an expensive high-dimensional integral.
There are dimension-in\-dependent MCMC methods \cite{hairer2014,vollmer2015dimension} that do not suffer from the curse of dimensionality.
However, these methods often require expensive pre-calculations.
Contrarily, naive proposal strategies may lead to computationally inefficient methods since the chain has to sequentially explore posterior probabilities in many dimensions.
In high dimensions one may need many evaluations of the likelihood (and thus forward model runs) in order to produce enough effective samples to adequately explore the posterior.
It can be infeasible to solve the model enough times to do so if it is computationally expensive.
The methane hydrate model is computationally expensive and has a somewhat high-dimension parameter space, so it is desirable to reduce the expense of MCMC by dimension reduction.
\subsection{Active subspaces}
\textit{Active subspaces} is a recently developed method for dimension reduction \cite{constantine2014computing,constantine2015active,constantine2014active} which identifies important directions in the parameter space.
It has been shown to be useful in several applications including approximation, integration, optimization, and sensitivity analysis \cite{constantine2015active,grey2017active,jefferson2016reprint,lukaczyk2014active}.
Recently, active subspaces have been used to reduce the dimension of parameter spaces in Bayesian inversion \cite{constantine2016accelerating,cortesi2017forwardbackward}.

To identify important directions of a function of interest $f:\R^\dim\to\R$, one looks at eigenpairs of the positive semi-definite $\dim\times\dim$ matrix 
\begin{align}
	\label{eq:C}
	\C &\defas \expct{\rho}{\grad f \grad f \tr} \\
	&= \int{\grad f(\x) \grad f(\x)\tr \rho(\x) \dx} = \bar{\W}\bar{\MeigVals}\bar{\W}\tr,
\end{align}
where $\rho$ is a given probability density function, and $f$ is continuous and differentiable on the support of $\rho$ and has square-integrable derivatives with respect to $\rho$.
In the Bayes\-ian inversion setting, $f$ is the data misfit function $\misfit$.
Since
\begin{equation}
	\bar{\eigVal}_i = \bar{\eigVec}_i\tr\C\bar{\eigVec}_i = \int{(\bar{\eigVec}_i\tr\grad f(\x))^2\rho(\x)\dx}, \quad i=1,\ldots,\dim,
\end{equation}
with $[\bar{\eigVec}_1,\ldots,\bar{\eigVec}_\dim]=\W$, it is evident that for small or even zero-valued eigenvalues the by $\rho$ weighted variation of $f$ is small (or even zero) in directions of corresponding eigenvectors, on average.
These directions can be treated differently than more relevant ones.

Matrix~$\C$ can be approximated by Monte Carlo integration:
\begin{equation}
	\label{eq:mc_sum}
	\C \approx \frac{1}{N}\sum_{j=1}^{N}{\grad f(\x_j)\grad f(\x_j)}\tr = \W\MeigVals\W\tr
\end{equation}
for samples $\x_j$ distributed according to $\rho$.
Matrices $\W$ and $\MeigVals$ denote perturbed versions of $\bar{\W}$ and $\bar{\MeigVals}$ due to the finite Monte Carlo sum in \eqref{eq:mc_sum}.
Since one actually is only interested in approximating the eigenpairs of $\C$ accurately, it is in most cases enough to have
\begin{equation}
	\label{eq:heur_N}
	N = \alpha\times\ell\times\log(\dim)
\end{equation}
samples, where $\alpha\in[2,10]$ is a \textit{sampling factor}, $\ell\leq\dim$ is the number of eigenpairs to be accurately approximated, and $\dim$ denotes the dimension of the function of interest $f$. For more details, we refer to \cite{constantine2014computing,constantine2015active}.
Gradients can be calculated by using methods such as adjoint approaches, finite differences and radial basis functions.

If there is a large enough spectral gap after the $k$th eigenvalue of $\C$, one can define a corresponding $k$-dimensional active subspace and variables on it.
Let $\W=(\W_1\;\;\W_2)$, where $\W_1\in\R^{\dim\times k}$ and $\W_2\in\R^{\dim\times \dim-k}$.
The first $k$ eigenvectors make up $\W_1$.
The active subspace is defined by the range of $\W_1$.
With enough samples in the Monte Carlo approximation of $\C$, a larger spectral gap leads to a more accurate approximation of the active subspace \cite{constantine2014computing}.
A more explicit bound on the number of samples $N$ has been recently published in \cite{holodnak2018probspacebound}.
The approximation quality of the active subspace can be estimated by using \textit{bootstrap intervals}.
We regard the quality of approximation in terms of so-called \textit{subspace estimation errors}.
The subspace estimation error is defined by
\begin{equation}
\text{dist}\left(\text{ran}(\bar{\W}_1),\text{ran}(\W_1)\right) = \norm{\bar{\W}_1\bar{\W}_1\tr - \W_1\W_1\tr}_2.
\end{equation}
For algorithmic details, i.e. how to compute the subspace errors, see \cite{constantine2014computing,constantine2015active}.

The input variable $\x$ can, according to the definition of an active subspace, be separated into an active and inactive part by
\begin{equation}
	\label{eq:W1yW2z}
	\x = \W\W\tr \x = \W_1\W_1\tr \x + \W_2\W_2\tr \x = \W_1\y + \W_2\z
\end{equation}
for $\y \defas \W_1\tr \x$ and $\z \defas \W_2\tr \x$.
The variable $\y\in\R^k$ is called the \textit{active}, and $\z\in\R^{\dim-k}$ is called the \textit{inactive variable}.

A lower-dimensional function $g:\R^k\to\R$, approximating $f$, can be constructed via a conditional expectation or a regression surface on the active variable \cite{constantine2015active}.
Because it would be very expensive to compute a conditional expectation in every MCMC step, we choose to find a regression function $g$ such that
\begin{equation}
	f(\x) \approx g(\W_1\tr \x) = g(\y)
\end{equation}
by using Algorithm~\ref{alg:regr_surf}.
\begin{algorithm}
	\caption{Computing the regression surface $g$ in the active variable}
	\label{alg:regr_surf}
	\begin{flushleft}
		Assume samples $\x_i$ $(i=1,\ldots,N)$ according to $\rho$ and corresponding function values $f_i$ $(i=1,\ldots,N)$ are given.
		\begin{enumerate}
			\item Compute samples $\y_i$ in the active subspace by
			\begin{equation}
				\y_i = \W_1\tr \x_i, \quad i=1,\ldots,N.
			\end{equation}
			\item Find a regression surface $g$ for pairs $(\y_i,f_i)$ such that
			\begin{equation}
				g(\y_i) \approx f_i, \quad i=1,\ldots,N.
			\end{equation}
			\item Get a low-dimensional approximation of $f$ by computing
			\begin{equation}
				f(\x) \approx g(\W_1\tr\x).
			\end{equation}
		\end{enumerate}
	\end{flushleft}
\end{algorithm}
\subsection{Active subspaces for MCMC}
\label{subsec:as_mcmc}
The low-dimensionality of $\misfit$ defined in Eqn.~\eqref{eq:misfit} can be exploited to accelerate MCMC.
For Bayesian inversion, the matrix $\C$ is computed with $f=\misfit$ and $\rho=\prior$.
The gradient of $\misfit$, which is needed for estimating $\C$, is given by
\begin{equation}
	\label{eq:grad_misfit}
	\grad\misfit(\x) = \grad\model(\x)\tr\cov\inv(\model(\x)-\data),
\end{equation}
where $\grad\model:\R^\dim\to\R^{\dimData\times n}$ is the Jacobian of the forward map $\model$.

Note that using the conditional expectation for approximating $f_\data$, it is possible to prove an upper bound on the Hellinger distance between the true and the corresponding approximating posterior \cite{constantine2016accelerating}.
More specific, defining
\begin{equation}
g_\data(\y) \defas \int_{\R^{n-k}}{f_\data(\W_1\y+\W_2\z)\,\rho_{\z|\y}(\z|\y)\dz},
\end{equation}
we have
\begin{equation}
\label{eq:bnd_hell}
d_H(\posterior,\rho_{\text{post},g_\data}) \leq C\left(\lambda_{k+1}+\cdots+\lambda_n\right)^{1/2},
\end{equation}
where
\begin{equation}
\rho_{\text{post},g_\data}(\x|\data) \defas \frac{\exp(-g_\data(\W_1\tr\x))\prior(\x)}{\int{\exp(-g_\data(\W_1\tr\x'))\prior(\x'){\dx'}}}
\end{equation}
and $C>0$.
That means small eigenvalues corresponding to the inactive subspace directly lead to a good approximation of the posterior (using the conditional expectation $g_\data$ as an approximation to $f_\data$).

The advantage of having a low-dimensional active subspace for the data misfit function $\misfit$ is that the Markov chain's mixing is accelerated when applying MCMC in few dimensions \cite{constantine2016accelerating}.
The Metropolis-Hastings algorithm can be adjusted to investigate only in the active subspace (see Algorithm~\ref{alg:as_mh_mcmc}), i.e. in the directions where the data misfit changes more, on average.
Equivalently, the Markov chain updates along the directions that are most informed by the data, on average.

\begin{algorithm}
	\caption{Metropolis-Hastings in the active variable $\y$}
	\label{alg:as_mh_mcmc}
	\begin{flushleft}
		Assume a symmetric proposal density function $\tau$, an initial point $\y_1$, a kernel density estimate $\hat{\rho}$ for the prior on the active variable $\priory$ and a response surface $g_\data$ for $\misfit$ in the active subspace are given.
		
		For $i=1,2,\ldots,N_\y-1$
		\begin{enumerate}
			\item Draw a proposal $\tilde{\y}$ from $\tau$ centered at $\y_i$.
			\item Calculate the acceptance ratio
			\begin{equation}
				\alpha(\tilde{\y},\y_i) = \min\left(1, \frac{\exp(-g_\data(\tilde{\y}))\hat{\rho}(\tilde{\y})}{\exp(-g_\data(\y_i))\hat{\rho}(\y_i)}\right).
			\end{equation}
			\item Draw $u\sim\mathcal{U}([0,1])$.
			\item Set $\y_{i+1} = \tilde{\y}$ if $\alpha(\tilde{\y},\y_i) \geq u$, otherwise set $\y_{i+1}=\y_i$.
		\end{enumerate}
	\end{flushleft}
\end{algorithm}

After applying Algorithm~\ref{alg:as_mh_mcmc}, one has samples $\y_i$ in the active subspace, i.e. \textit{active samples}.
These samples are naturally correlated due to the Markov property of the chain.
The minimum effective sample size $\NyEss\defas\min_{l\in\lbrace1,\ldots,k\rbrace}\NylEss$ can be found by taking the minimum of all effective sample sizes in the components of samples $\y_i$.
The effective sample size $\NylEss$ for the $l$th component in samples $\y_i$ is computed via the formula \cite{brooks2011handbook}
\begin{equation}
	\label{eq:ess}
	\NylEss = \frac{N_\y}{1+2\sum_{j=1}^{\Jmax}{r^{(l)}_j}},
\end{equation}
where $N_\y$ is the total amount of samples in the chain (the chain length), $\Jmax$ stands for the maximum lag taken into account, and $r^{(l)}_j$ denotes the correlation between the $l$th component of samples $\y_i$ with lag $j$.
As it is common in time series analysis and digital signal processing, the correlations $r^{(l)}_j$ are calculated using (inverse) fast Fourier transforms.
Take $\NyEss$ out of $N_\y$ samples with equally distance and regard them as the (effective) result of sampling from the posterior in the active subspace.

In order to get (effective) samples $\x_i$ in the original space, it is necessary to also produce samples in the inactive subspace, conditioned on the respective active samples.
The conditional density $\priorzy(\z|\y)$ of $\z$ given $\y$ is
\begin{align}
	\label{eq:prizy}
	\begin{split}
		\priorzy(\z|\y) ={}& \frac{\prior(\W_1\y+\W_2\z)}{\int{\prior(\W_1\y+\W_2\z')\dz'}} \\
		\propto{}& \prior(\W_1\y+\W_2\z).
	\end{split}
\end{align}
Hence, for every effective active sample $\y_i$, we additionally compute $N_{\z,i}$ inactive samples $\z^{(i)}_j$ $(j=1,\ldots,N_{\z,i})$ by applying Algorithm~\ref{alg:ias_mh_mcmc} to Eqn.~\eqref{eq:prizy}.
\begin{algorithm}
	\caption{Metropolis-Hastings in the inactive variable conditioned on active sample $\y$}
	\label{alg:ias_mh_mcmc}
	\begin{flushleft}
		Assume an active sample $\y$, a symmetric proposal density function $\tau$, and an initial point $\z_1$ is given.

		For $j=1,2,\ldots,N_\z-1$
		\begin{enumerate}
			\item Draw a proposal $\tilde{\z}$ from $\tau$ centered at $\z_j$.
			\item Calculate the acceptance ratio
			\begin{equation}
				\alpha(\tilde{\z}, \z_j) = \min\left(1, \frac{\prior(\W_1\y+\W_2\tilde{\z})}{\prior(\W_1\y+\W_2\z_j)}\right).
			\end{equation}
		\item Draw $u\sim\mathcal{U}([0,1])$.
		\item Set $\z_{j+1} = \tilde{\z}$ if $\alpha(\tilde{\z},\z_j) \geq u$, otherwise set $\z_{j+1}=\z_j$.
		\end{enumerate}
	\end{flushleft}
\end{algorithm}
Afterwards, we again compute the minimum effective sample size $\NziEss$ of samples~$\z^{(i)}_j$ similar to Eqn.~\eqref{eq:ess}.
Note that the effective sample sizes $\NziEss$ can be different.
So, in order to get an equal number $\NzEss$ of effective inactive samples $\z^{(i)}_j$ per effective active sample $\y_i$, we have to choose the chain lengths $N_{\z,i}$ appropriately (according to the corresponding autocorrelations).
By using Eqn.~\eqref{eq:W1yW2z}, we then get $\NxEss=\NyEss\times\NzEss$ effective samples $\x_i$ in the original parameter space.
These samples approximate the posterior density.

\section{Results}
\label{sec:results}
\begin{figure}
	\centering
	\includegraphics[scale=0.5]{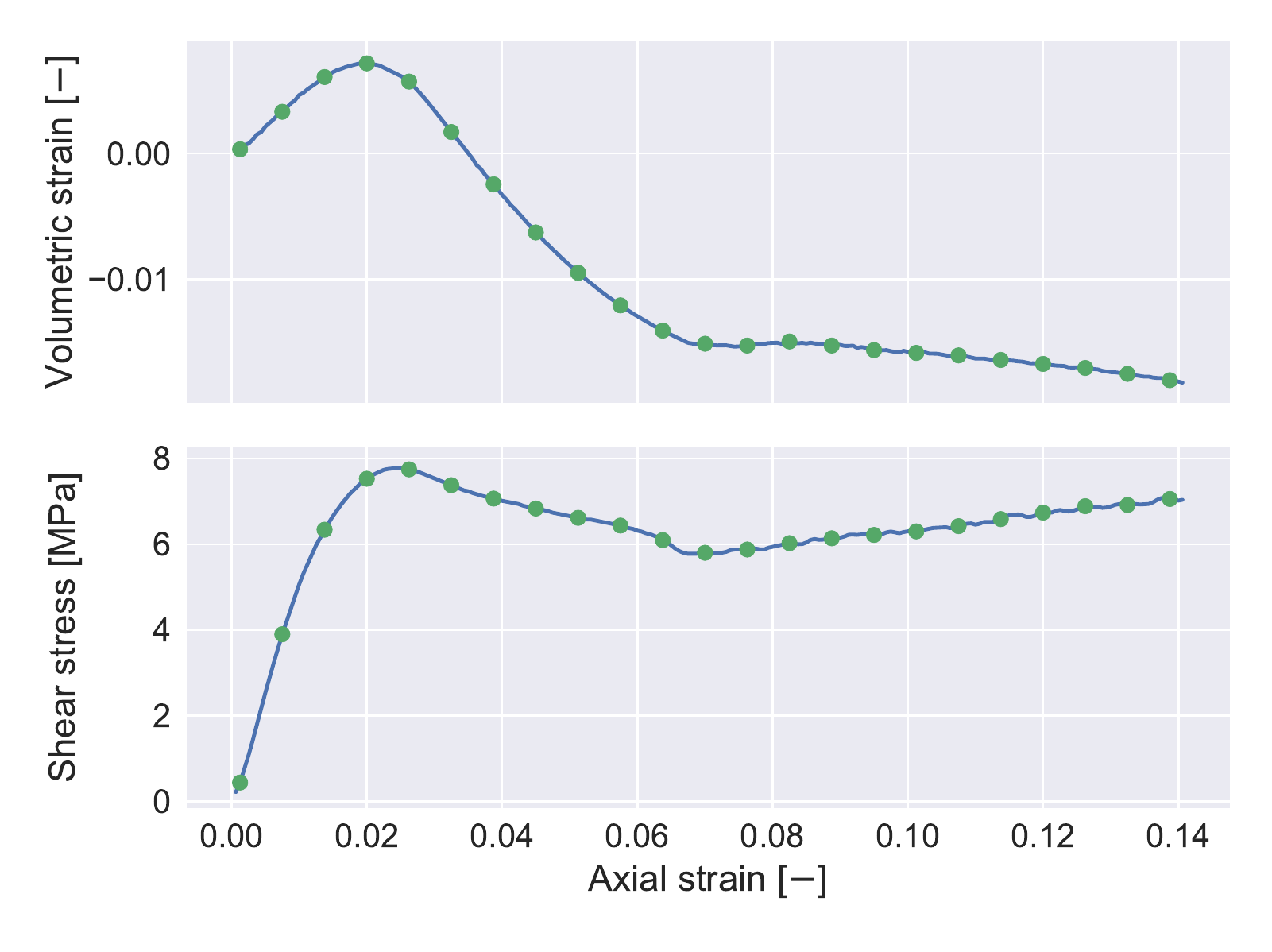}
	\caption{Volumetric strain $\strain$ and shear stress $\stress$ data plotted versus axial strain.}
	\label{fig:qoi}
\end{figure}

The efficient inversion method from Section~\ref{sec:inversion} is used to infer parameters of the model described in Section~\ref{sec:model}.
The quantities of interest are the volumetric strain, $\strain$, and shear stress, $\stress$, for 23 given axial strain values.
Fig.~\ref{fig:qoi} shows the measured values representing the data, $\data$.
The forward map is 
\begin{equation}
	\label{eq:fwd_model}
	\model(\x) = \begin{pmatrix}\modelStrain(\x) \\ \modelStress(\x)\end{pmatrix}\in\R^{\dimStrain+\dimStress},
\end{equation}
where $\dimStrain=\dimStress=23$, and values of $\modelStrain:\R^8\to\R^\dimStrain$ and $\modelStress:\R^8\to\R^\dimStress$ are ordered according to the corresponding axial strain.
Note that $\dim=8$ and $\dimData=\dimStrain+\dimStress$ in the notation of Section~\ref{sec:inversion}.
The data is
\begin{equation}
	\data=\begin{pmatrix}\dataStrain \\ \dataStress\end{pmatrix}\in\R^{\dimStrain+\dimStress}.
\end{equation}
The covariance of the Gaussian noise~$\noise$ has the form
\begin{equation}
	\cov = \begin{pmatrix}\covStrain & 0 \\ 0 & \covStress\end{pmatrix},
\end{equation}
where $\covStrain$ and $\covStress$ are diagonal matrices corresponding to a 2\% noise level for each measurement.
The data misfit function is
\begin{equation}
	\label{eq:misfits}
	\misfit(\x) = \frac{\norm{\dataStrain-\modelStrain(\x)}_\covStrain^2}{2} + \frac{\norm{\dataStress-\modelStress(\x)}_\covStress^2}{2}.
\end{equation}

\subsection{Computational costs}
To calculate the active subspace, $\C$ must be approximated with Eqn.~\eqref{eq:mc_sum}.
In Eqn.~\eqref{eq:grad_misfit}, the Jacobian of the forward forward map $\model$ is needed to calculate the gradient of the data misfit.
Since the model is not readily capable of adjoint computations, central finite differences are used to approximate partial derivatives
\begin{equation}
	\frac{\partial \model}{\partial \x^{(j)}} \approx \frac{\model(\x^{(j)}+h \mathbf{e}_j)-\model(\x^{(j)}-h \mathbf{e}_j)}{2h}
\end{equation}
for $j=1,\ldots,8$.
For one sample of $\grad\misfit$, forward evaluations of $\model$ must be computed $2\times 8+1=17$ times.
Using the heuristic from Eqn.~\eqref{eq:heur_N}, one needs at least $N=\ceil{\alpha\times8\times\log(8)} = 167$ gradient samples to estimate $\C$ (with the most pessimistic choice of $\alpha=10$).
In order to compute the subspace estimation errors more accurately, we use $N=250$ samples.
Fortunately, the calculation of the gradient samples is trivially parallelizable.
To compute $17\times 250$ samples, we needed $15.12$ hours on $35$ cores, i.e. $529.2$ core hours.
One forward evaluation therefore took $7.47$ minutes on average.
The additional calculations required to estimate the active subspace are very cheap.
\subsection{Active subspace for the model}
We assume a uniform prior distribution on the hypercube $[-1,1]^8$, i.e. $\prior=\mathcal{U}([-1,1]^8)$. The samples are then mapped to the intervals in Table~\ref{tab:prior_intervals} before running the forward model.
It is common to calculate the active subspace with respect to a mean-zero, scaled prior.
The prior intervals in Table~\ref{tab:prior_intervals} were chosen using sensitivity analysis and engineering knowledge of the model.
Sensitivity analysis indicates that there is a promising local minimum of the data misfit function in these intervals.
The prior intervals for parameters $\lambdadotmax$ and $\lambdamax$ are chosen to be small since the model is extremely sensitive to those parameters, and values outside of these ranges greatly diverge from the data.

\begin{table}
	\caption{Prior intervals of the eight parameters.}
	\label{tab:prior_intervals}
	\centering
	\begin{tabular}{lllll}
		\hline\noalign{\smallskip}
		No. & Parameter & Min & Max & Unit \\
		\hline\noalign{\smallskip}
		1 & $c$              & $1.8\times10^6$ & $2.4\times10^6$ & Pa \\
		2 & $\alpharesmin$   & $0.5$ & $0.6$ & -- \\
		3 & $\deltaalphares$ & $0.2$ & $0.3$ & -- \\
		4 & $\lambdadotmax$  & $1.6\times10^{-3}$ & $1.9\times10^{-3}$ & -- \\
		5 & $\malpha$        & $0.75$ & $1.05$ & -- \\
		6 & $\betamax$       & $0.3$ & $0.45$ & -- \\
		7 & $\lambdamax$     & $0.01$ & $0.011$ & -- \\
		8 & $\mbeta$         & $0.67$ & $0.74$ & -- \\
		\hline\noalign{\smallskip}
	\end{tabular}
\end{table}

\begin{figure}
	\centering
	\includegraphics[scale=0.5]{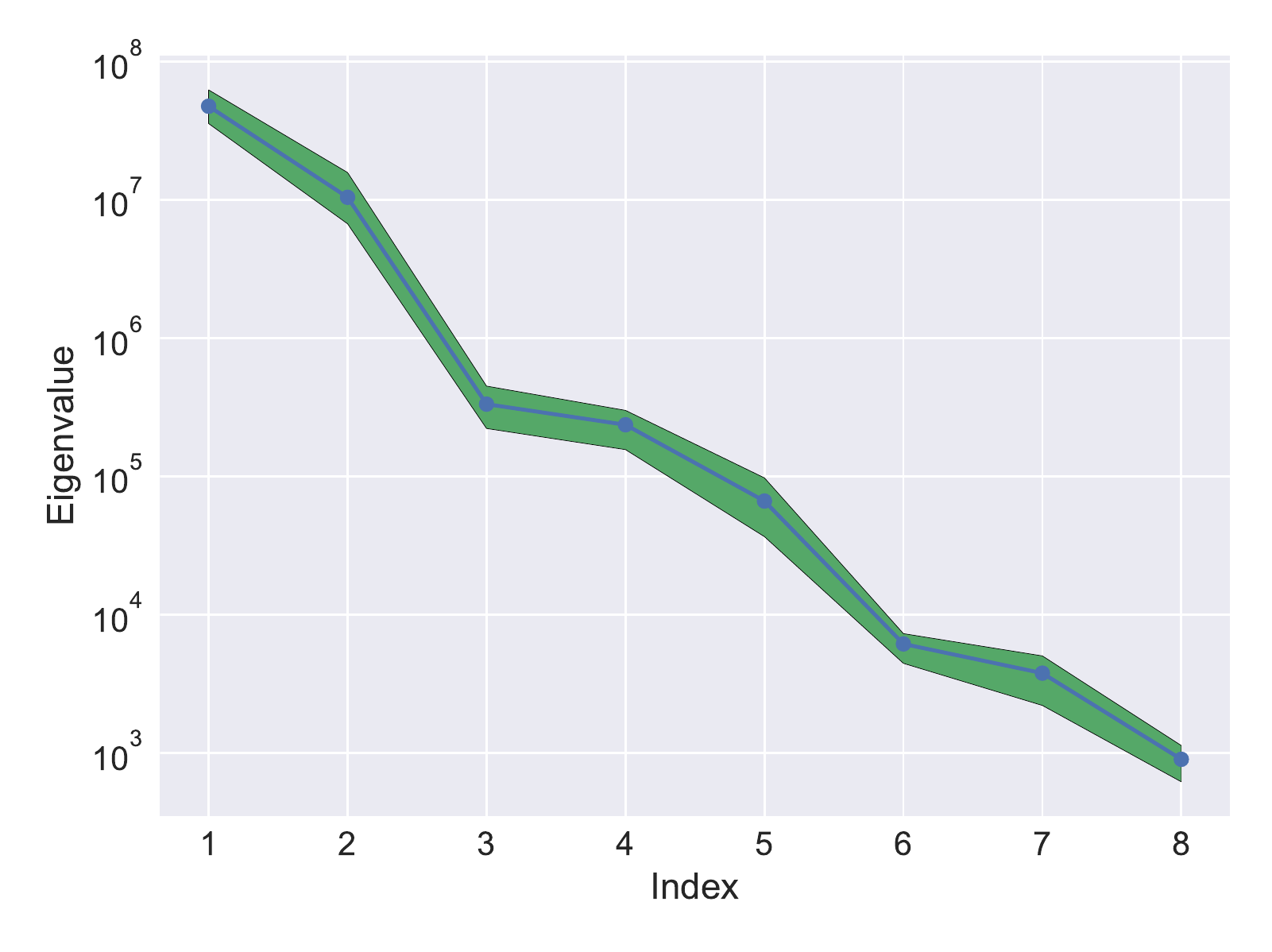}
	\caption{Approximated eigenvalues of $\C$ with bootstrap intervals.}
	\label{fig:eigvals}
\end{figure}

\begin{figure}
	\centering
	\includegraphics[scale=0.5]{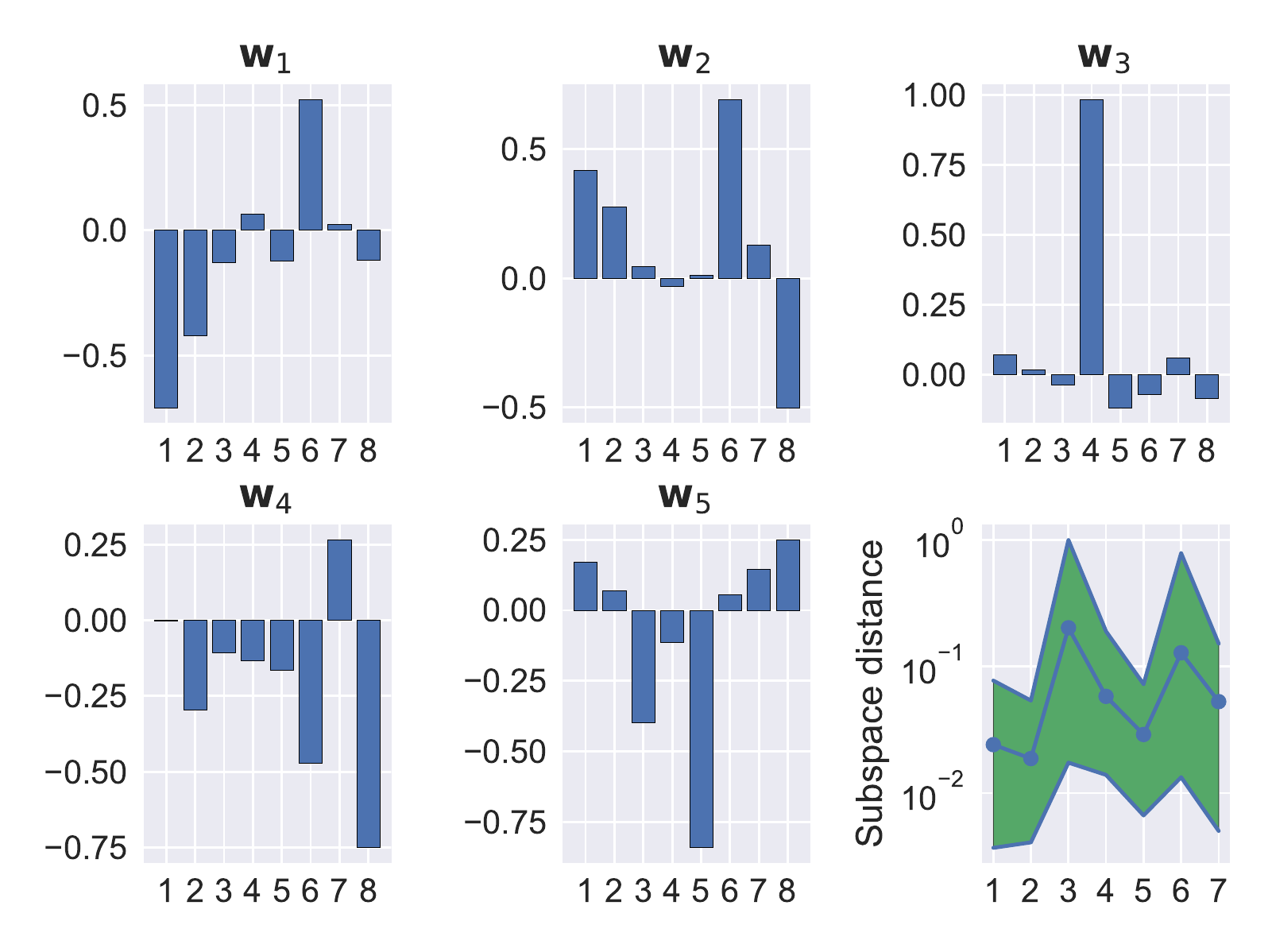}
	\caption{Components of the first five eigenvectors of $\C$.
		The lower right plot shows the subspace distance for each subspace.}
	\label{fig:eigvecs}
\end{figure}

The eigenvalue decomposition of the approximation of $\C$ using the 250 data misfit gradient samples results in the eigenvalues shown in Fig.~\ref{fig:eigvals}.
Fig \ref{fig:eigvecs} shows the components of the first five eigenvectors and the subspace errors for each subspace from using the Monte Carlo sum from Eqn.~\eqref{eq:mc_sum} \cite{constantine2014computing,constantine2015active}.
To identify active subspaces, we look for gaps in the spectrum.

We see that rather large gaps appear after the second and fifth eigenvalues.
This fact is confirmed by the corresponding subspace estimation errors (lower right plot in Fig.~\ref{fig:eigvecs}) that are very small for 2D and 5D spaces.
The eigenvectors associated with the larger eigenvalues define the directions in the parameter space that are more informed by data, on average.

\begin{figure}
	\centering
	\includegraphics[scale=0.5]{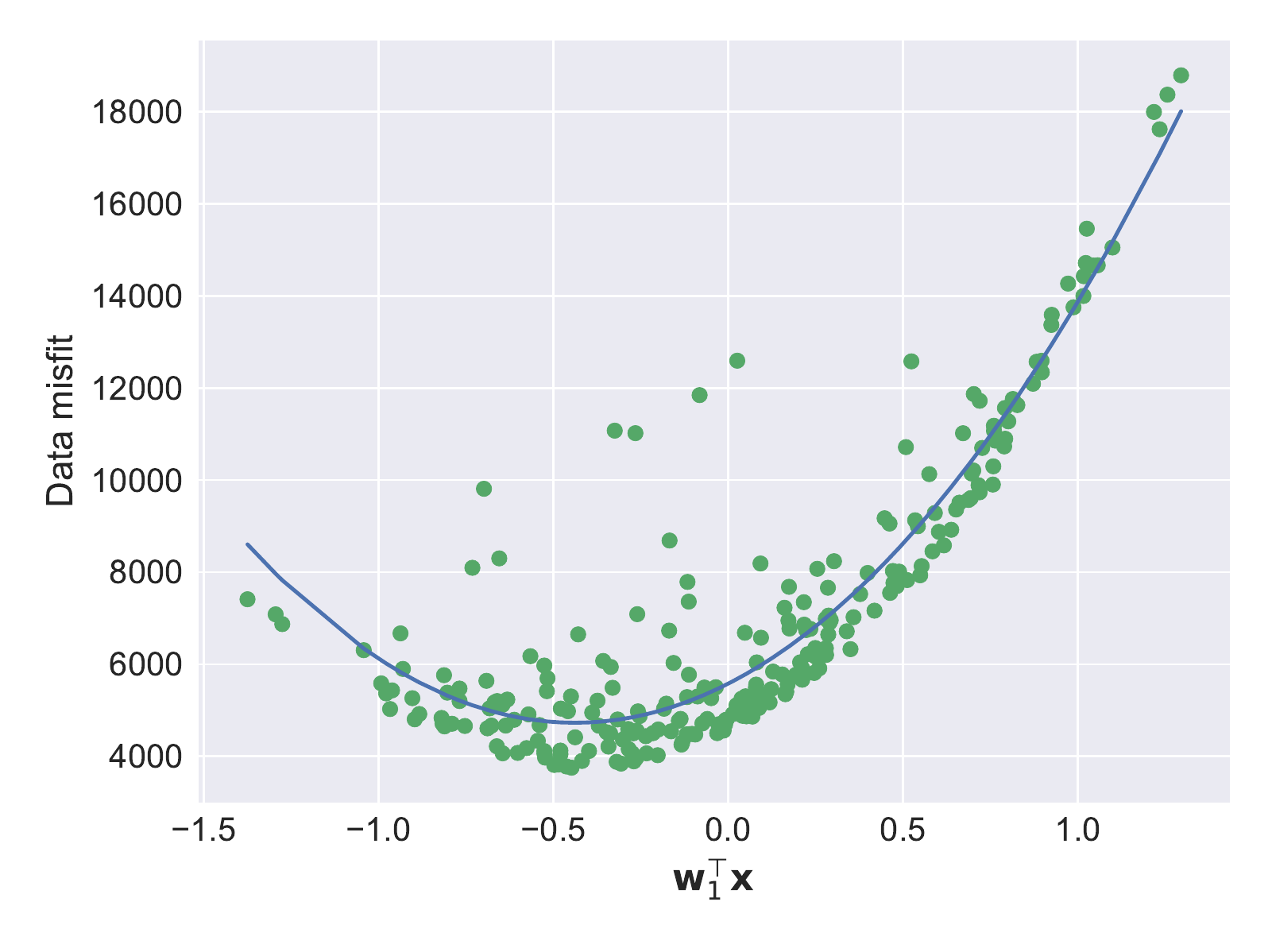}
	\caption{1D summary plot of $\misfit$ on the axis $\eigVec_1^\top\x$. The quadratic regression fit has a coefficient of determination of $r^2=0.8434$.}
	\label{fig:sum_plot_1d}
\end{figure}

\begin{figure}
	\centering
	\includegraphics[scale=0.5]{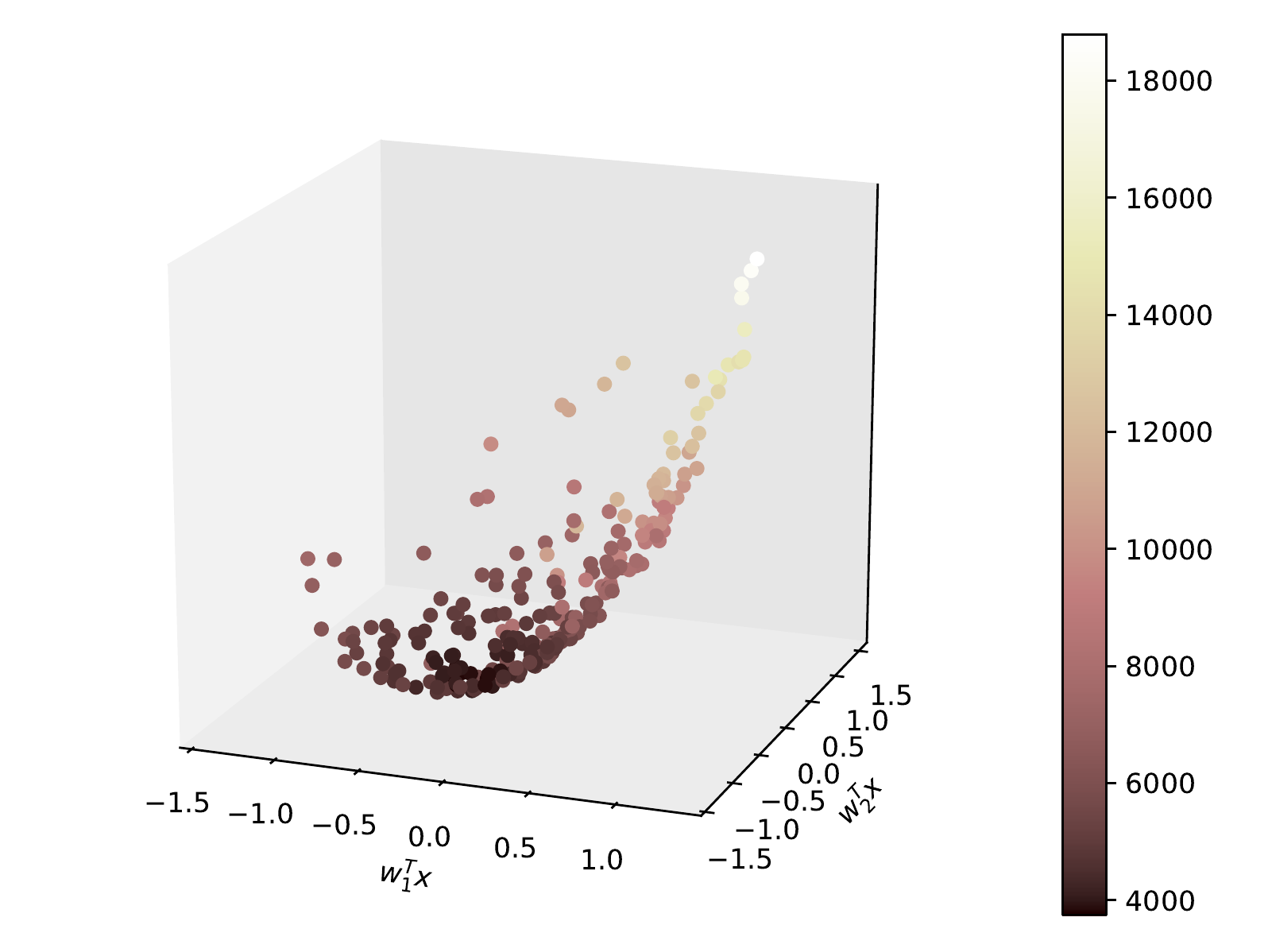}
	\caption{2D summary plot of $\misfit$ on the axes $\eigVec_1^\top\x$ and $\eigVec_2^\top\x$. The quadratic regression fit has a coefficient of determination of $r^2=0.9776$.}
	\label{fig:sum_plot_2d}
\end{figure}

1D and 2D summary plots for the first two eigenvectors are plotted in Figs.~\ref{fig:sum_plot_1d} and \ref{fig:sum_plot_2d}.
They show the data misfit function on the axes $\eigVec_1\tr\x$ and $\eigVec_2\tr\x$, i.e. the coordinate system is rotated and projected onto the subspace to see if there exists a lower-dimensional structure.
The response surfaces are constructed according to Algorithm~\ref{alg:regr_surf}.
A quadratic polynomial was fitted to the data misfit values in the active subspace with the Python package \textit{SciKit Learn} \cite{scikit-learn}.
The 1D summary plot already indicates low-dimensionality but contains some outliers resulting in a coefficient of determination (or $r^2$ score)  of $r^2=0.8434$ which is not ideal.
The 2D summary plot shows a strong two-dimensional structure as is confirmed quantitatively by a coefficient of determination of $r^2=0.9776$.
This verifies what was indicated by the large gap after the second eigenvalue: that there is a strong 2D active subspace. 
For comparison, using a 5D active subspace (corresponding to the second large gap in the spectrum) leads to $r^2=0.9824$.
Both the 2D and the 5D regressions are sufficient for further use for MCMC.

Such regression-based surrogates are extremely cheap to evaluate, even in higher dimensions.
There are advantages to both the higher and lower-dimensional active subspaces.
The regression-based surrogate for the higher-dimensional space introduces less error into the system; however, a lower-dimensional active subspace corresponds to accelerated mixing for MCMC \cite{constantine2016accelerating}.

\subsection{Inference via MCMC in the active subspace}
After the low-dimensional representations of the data misfit function are calculated, they can be exploited within MCMC as was discussed Section~\ref{subsec:as_mcmc}.
Since there are two large gaps in the approximated eigenvalues of $\C$, we compare inference results for both the 2D and 5D active subspaces.

We expect that results using the 5D subspace approximate the exact posterior more accurately than the 2D subspace due to the bound on the Hellinger distance in Eqn.~\eqref{eq:bnd_hell}.
For the 5D subspace, fewer eigenvalues from the inactive subspace contribute to the upper bound.

In order to sample from the posterior distribution with MCMC, one not only needs a low-dimensional approximation of $\misfit$ but also a sufficiently accurate estimator $\prioryEst$ for the prior probability density $\priory$ marginalized on the active variable.
Since we have a multivariate uniform prior density, for which it is in general difficult to find analytical expressions or marginal densities, especially in higher dimensions, $\priory$ has to be estimated numerically.
We use a kernel density estimation (KDE) approach from SciKit Learn \cite{scikit-learn}.

A few trial runs of several thousand steps with different proposal covariances are used to determine ones which result in reasonable acceptance rates for the Metropolis-Hastings algorithm. 
The additional computational costs for these runs are small since the surrogate is very cheap.
There exist several more complicated techniques for finding appropriate proposal densities, but the cheapness of evaluating the surrogate allows for such trial and error, with each trial only taking a few minutes.
For MCMC in the 2D active subspace, we use a Gaussian proposal density function with covariance matrix $0.02\times \I$, and the acceptance rate is $35\%$.
The chain consists of $10^6$ steps with a burn-in of the first $10^5$ samples.
Using Eqn.~\eqref{eq:ess}, this results in an effective sample size of $\NyEss=75,000$ active posterior samples.
For the 5D case, we use a Gaussian proposal density with covariance matrix $0.0017\times \I$ in the active subspace and get an acceptance rate of again $35\%$ with $10^7$ steps and a burn-in of $10^6$ samples resulting in $21,000$ effective active samples.

Note that we do not directly exploit the quadratic nature of the data misfit function to further accelerate MCMC.
Such quadratic structure could potentially be used to compute proposal covariance matrices that accelerate mixing.
However, since this paper is intended as a showcase of active subspaces for Bayesian inference in geophysical applications, we do not further utilize the specific quadratic structure of the data misfit
Also, other setups and models might give greatly different shapes for the data misfit function.

\begin{figure}
	\centering
	\includegraphics[scale=0.5]{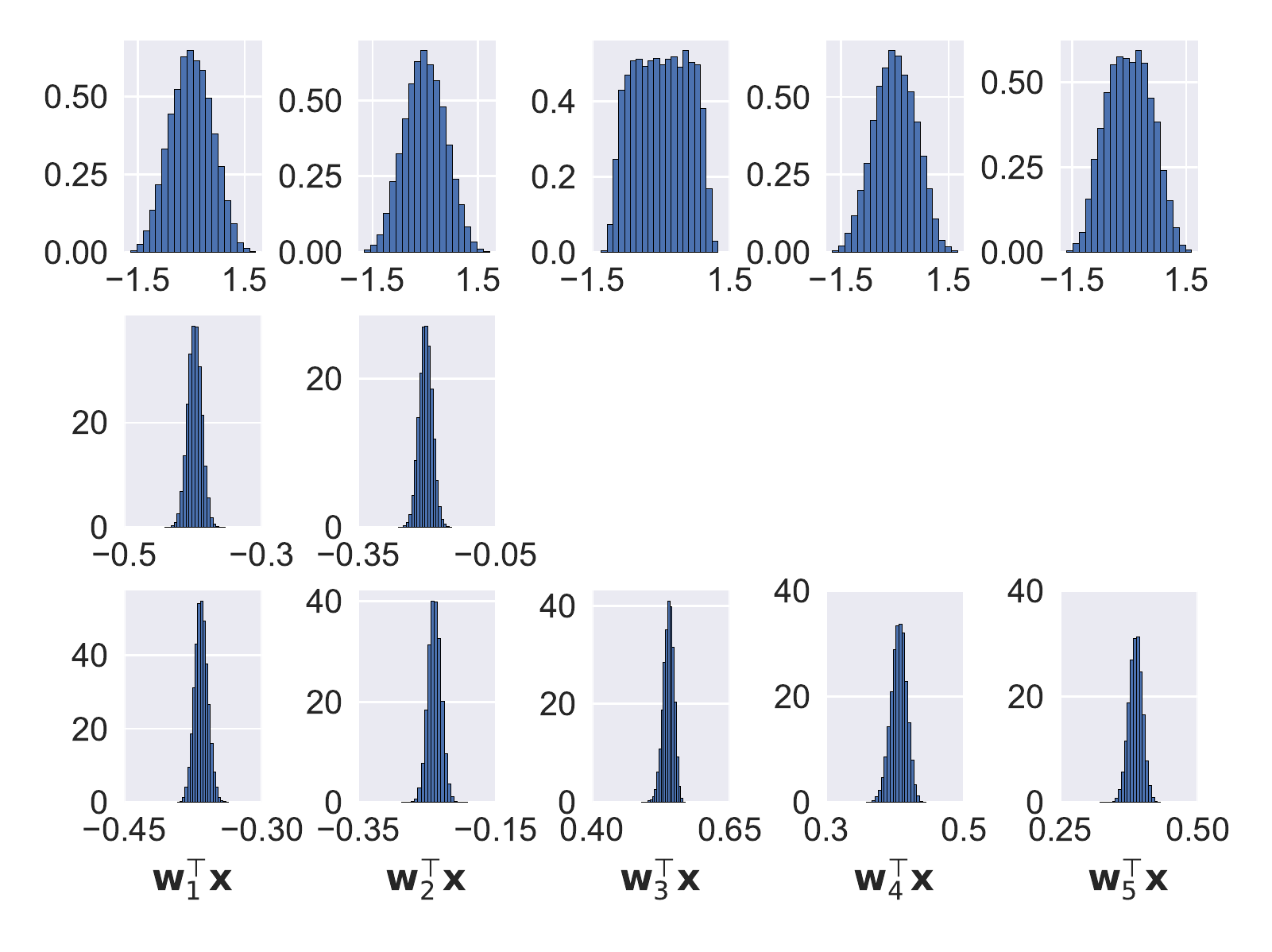}
	\caption{1D marginals of the prior (top row) and 2D (middle row) and 5D (bottom row) posteriors on the active variables.
		Note the different scales for respective $x-$ and $y-$axes.}
	\label{fig:prior_actvar}
\end{figure}

Fig.~\ref{fig:prior_actvar} shows histograms for the estimated prior density $\prioryEst$ in five dimensions and the posterior densities for the 2D and 5D subspaces.
Note that the two upper left plots show $\prioryEst$ in the 2D active subspace.
Also note that the scales of the $x-$ and $y-$axes differ in the subplots.
The different scales show a large reduction in variance from prior to posterior marginals.
That means that the data misfit was very informed in the directions of the active subspaces.
The active samples are concentrated in very small regions of the active subspaces.
Also, it is worthwhile to take a closer look to the first row displaying marginal prior densities in the active directions.
The histogram in the middle shows the prior marginalized for the 3rd active variable $\eigVec_3\tr\x$ which is approximately $\x_4$ as can be seen from Fig.~\ref{fig:eigvecs}.
This explains why the density still looks rather uniform compared to the other marginals in which more than one component shows up.
It suggests that active variables can not only be viewed as projections but also as a linear combination or weighted sum.
The prior densities certainly contain characteristics known from weighted sums of uniform random variables.

\begin{figure}
	\centering
	\includegraphics[scale=0.6]{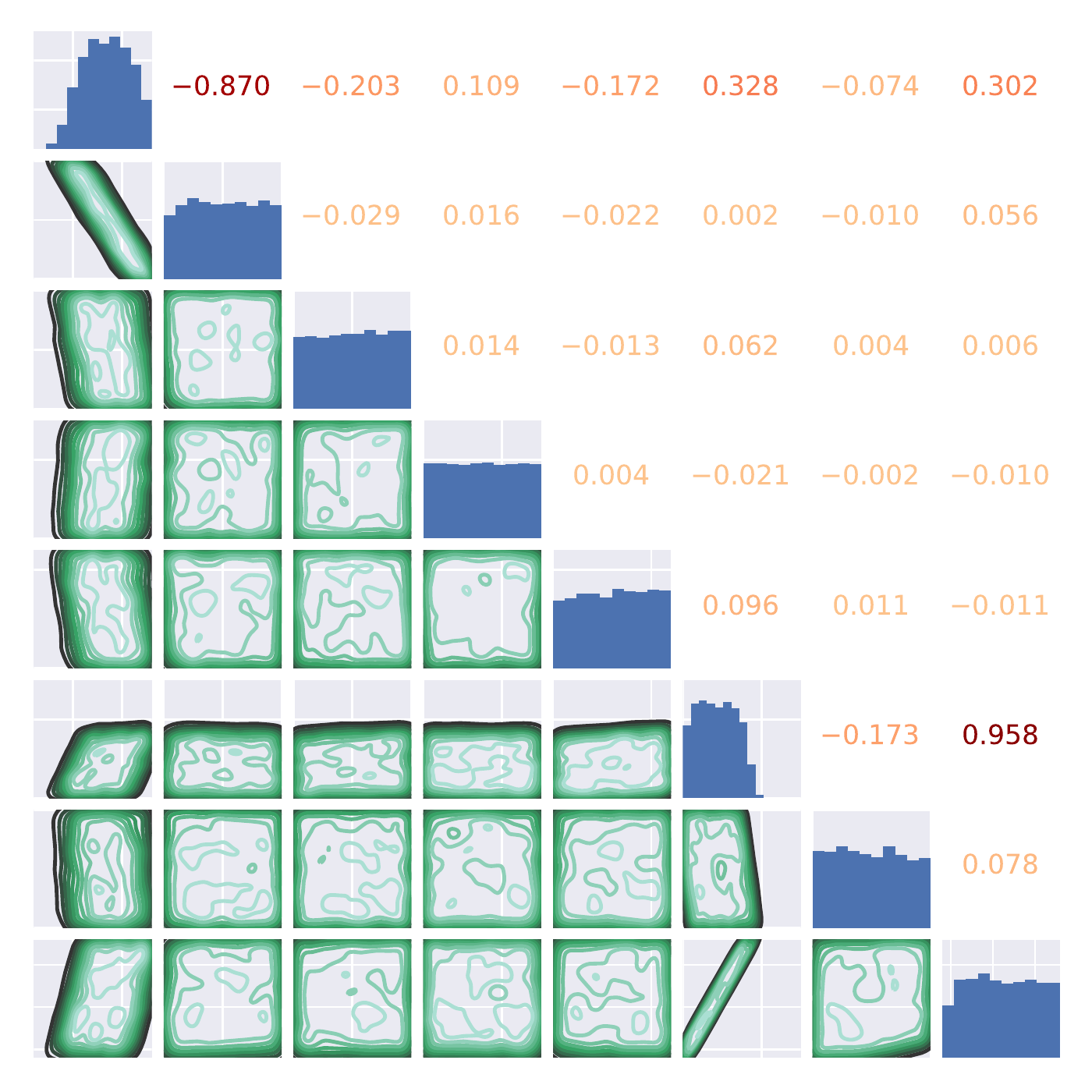}
	\caption{Marginal posterior densities for each parameter in the original space by using a 2D active subspace.
		The correlations in the upper triangle are colored according to their absolute value.}
	\label{fig:posterior_2d}
\end{figure}

\begin{figure}
	\centering
	\includegraphics[scale=0.6]{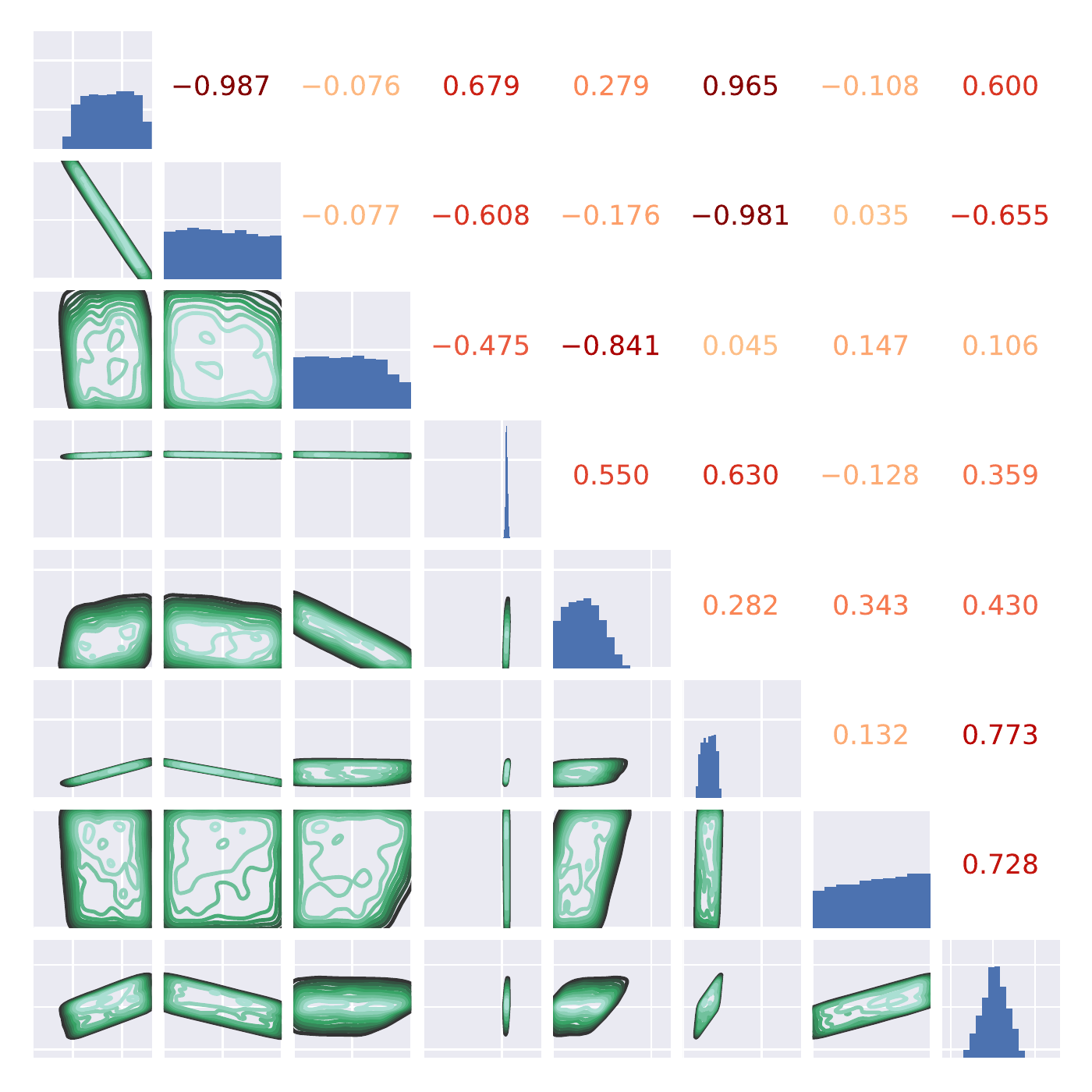}
	\caption{Marginal posterior densities for each parameter in the original space by using a 5D active subspace.
		The correlations in the upper triangle are colored according to their absolute value.}
	\label{fig:posterior_5d}
\end{figure}

\begin{figure}
	\centering
	\includegraphics[scale=0.5]{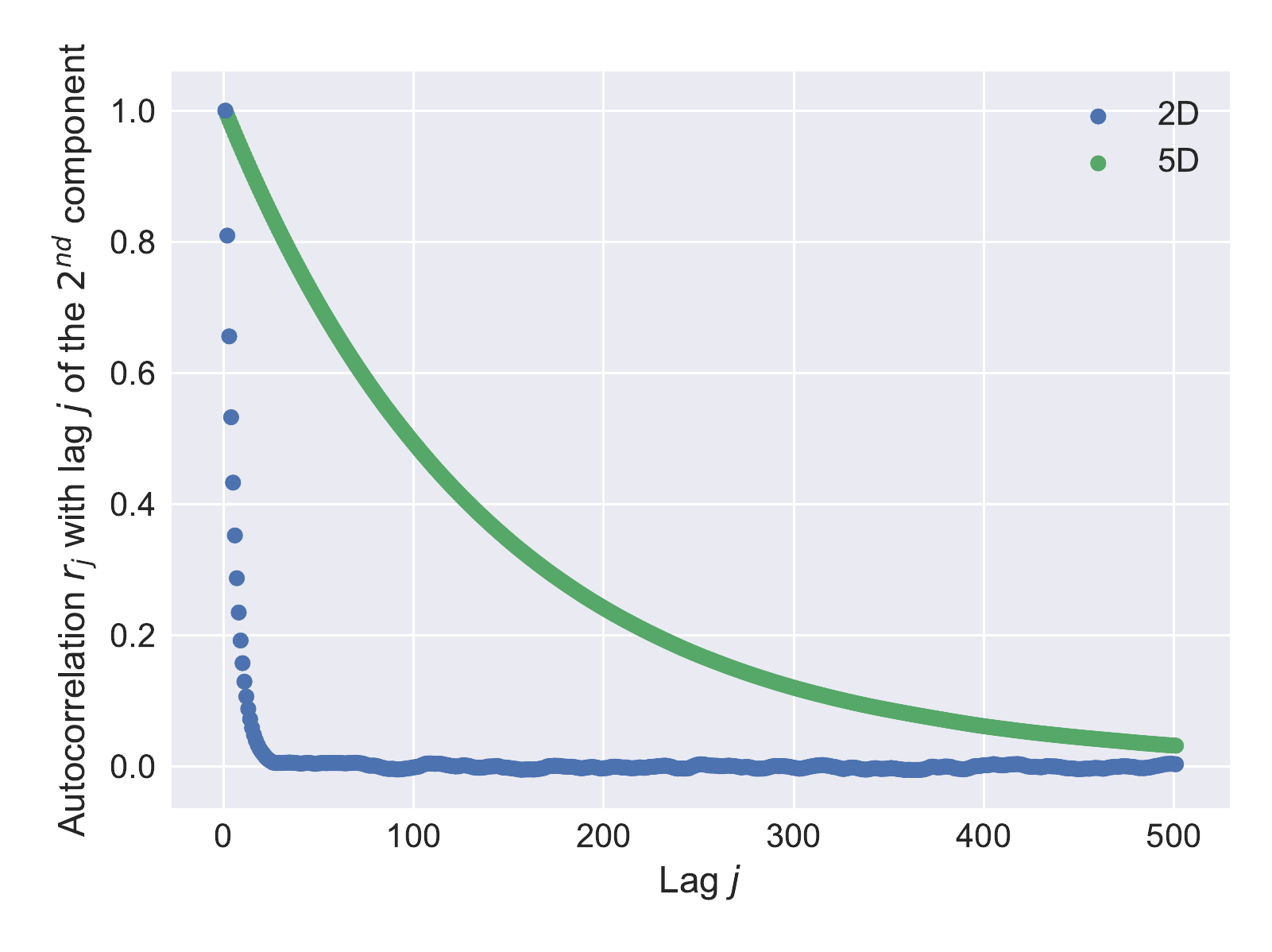}
	\caption{Autocorrelations in the second component of the MCMC chains for the active variable with a 2D and 5D active subspace.}
	\label{fig:autocorrs}
\end{figure}

\begin{figure}
	\centering
	\includegraphics[scale=0.5]{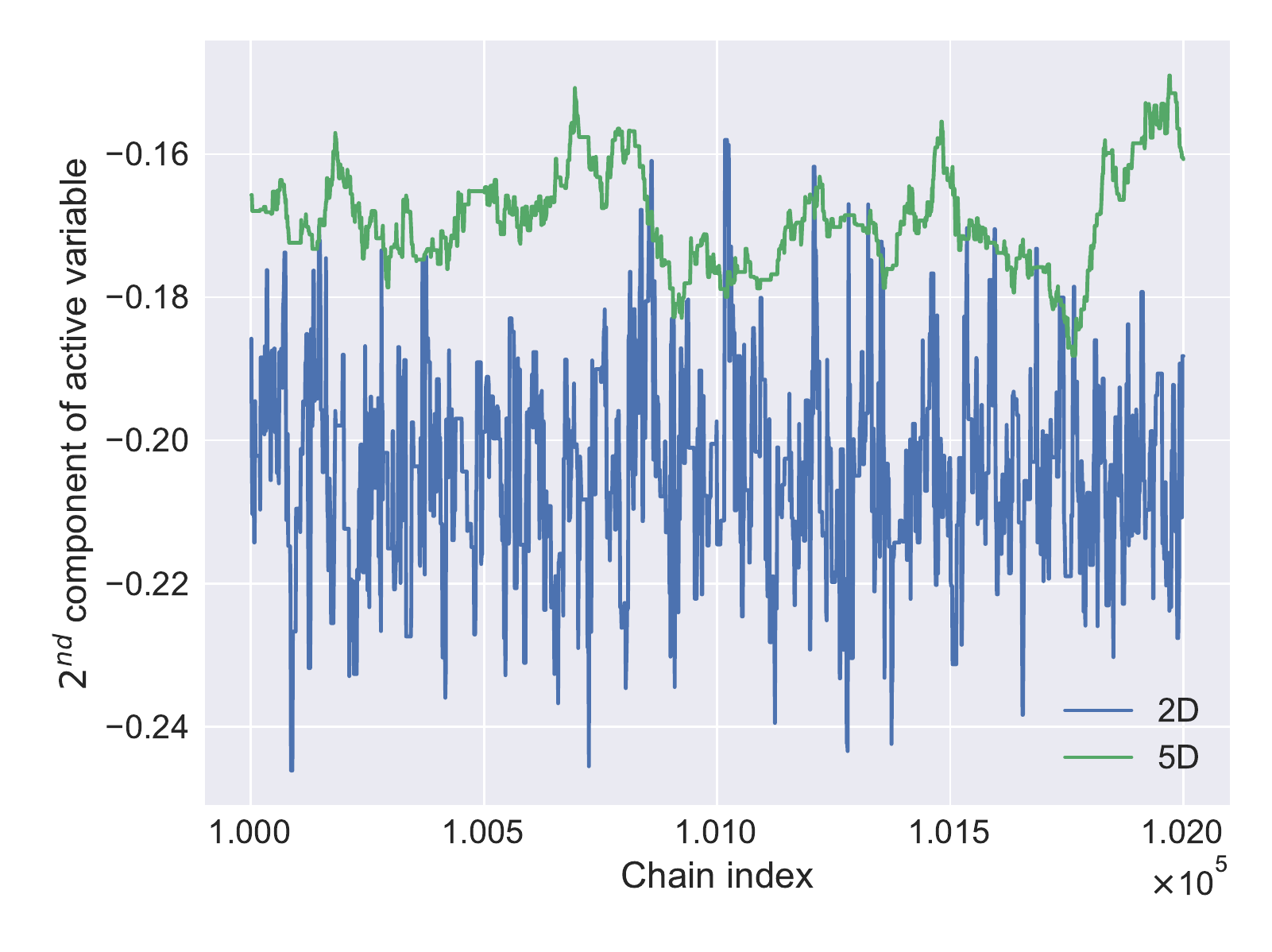}
	\caption{Mixing for the second component in the MCMC chains for the active variable with a 2D and 5D active subspace.}
	\label{fig:mixing}
	\end{figure}

We want to calculate samples of the posterior on the full parameter space.
Using Algorithm~\ref{alg:ias_mh_mcmc} for every effective active sample $\y_i$ and the decomposition from Eqn.~\eqref{eq:W1yW2z} gives effective samples $\x_i$.
For the 2D active subspace, we compute effective inactive samples $\z_i$ again using a Gaussian proposal density with covariance matrix $0.4\times \I$ with an acceptance rate of approximately $37\%$ depending on the  active sample.
The resulting marginal posterior densities are shown in Fig.~\ref{fig:posterior_2d}.
In the 5D subspace, a Gaussian proposal density with covariance matrix $0.8\times\I$ for computing effective inactive samples gives an acceptance rate of $36\%$.
The resulting marginal posterior densities are shown in Fig.~\ref{fig:posterior_5d}.
The resulting samples from MCMC computations are used to compute the posterior means and standard deviations of each parameter (see Table~\ref{tab:post_stats}).
The autocorrelations and mixing behavior for one component of the MCMC chains are plotted in Figs.~\ref{fig:autocorrs} and \ref{fig:mixing}.
The autocorrelations for the 5D subspace decrease slowly, meaning that the distance between two nearly uncorrelated samples in the chain is large.
For 2D, the situation is different since the autocorrelations are dropping rapidly.
These behaviors can be explained by the mixing within the chains.
The chains for the 5D space do not mix very well, which can be seen by rather small amplitudes in the oscillations.
Succeeding samples are quite correlated, resulting in longer chains in order to get enough uncorrelated (effective) samples.
In contrast, the samples from the 2D chain oscillate much giving more effective samples as in 5D (with the same number of samples).
Table~\ref{tab:perf_sampl} shows the quantitative performance of creating active samples in two and five dimensions, confirming what was mentioned earlier.
In two dimensions, we get more effective active samples per second.
However, since one step in the Markov chain is very cheap due to the response surface approximation, we are able to and choose to run the chain in five dimensions to get more informed posterior samples improving the inference result.
The fact that the chain has to be much longer is acceptable because of the cheap surrogate model.
The purpose of presenting the comparison between two dimensions is to once more point out the well-known potential behavior of Markov chains in different dimensions.

\begin{table}
	\centering
	\caption{Computational performance of creating active samples in different dimensions.}
	\label{tab:perf_sampl}
	\begin{tabular}{ccc}
		\hline\noalign{\smallskip}
		& 2D & 5D \\
		\hline\noalign{\smallskip}
		Eff. samples $\y_i$ & 90,910 & 3,301 \\
		\hline\noalign{\smallskip}
		Elapsed sec. & 1249.2 & 935.7 \\
		\hline\hline\noalign{\smallskip}
		\textbf{Eff. samples / sec.} & $\approx \mathbf{72.77}$ & $\approx \mathbf{3.53}$ \\
		\hline
	\end{tabular}
\end{table}

\begin{table*}
	\caption{Posterior means ($\vec{\mu}_{\text{post}}$) and standard deviations ($\sigma_\text{post}$) for each parameter using the 2D and 5D subspaces.}
	\label{tab:post_stats}
	\centering
	\resizebox{\linewidth}{!}{%
	\begin{tabular}{cccccccccc}
		\hline\noalign{\smallskip}
		& & $c$ & $\alpharesmin$ & $\deltaalphares$ & $\lambdadotmax$ & $\malpha$ & $\betamax$ & $\lambdamax$ & $\mbeta$ \\
		\noalign{\smallskip}\hline\noalign{\smallskip}
		\multirow{2}{*}{$2$D} & $\vec{\mu}_{\text{post}}$ & $2.169\times 10^6$ & $0.551$ & $0.251$ & $1.749\times 10^{-3}$ & $0.903$ & $0.343$ & $1.049\times 10^{-2}$ & $0.706$ \\
		& $\sigma_{\text{post}}$ & $1.165\times 10^5$ & $2.850\times 10^{-2}$ & $2.882\times 10^{-2}$ & $8.668\times 10^{-5}$ & $8.639\times 10^{-2}$ & $2.427\times 10^{-2}$ & $2.864\times 10^{-4}$ & $1.946\times 10^{-2}$ \\
		\noalign{\smallskip}\hline\noalign{\smallskip}
		\multirow{2}{*}{$5$D} & $\vec{\mu}_{\text{post}}$ & $2.183\times 10^6$ & $0.548$ & $0.246$ & $1.812\times 10^{-3}$ & $0.827$ & $0.332$ & $1.053\times 10^{-2}$ & $0.701$ \\
		& $\sigma_{\text{post}}$ & $1.148\times 10^5$ & $2.874\times 10^{-2}$ & $2.724\times 10^{-2}$ & $2.626\times 10^{-6}$ & $4.361\times 10^{-2}$ & $7.745\times 10^{-3}$ & $2.830\times 10^{-4}$ & $7.248\times 10^{-3}$ \\
		\noalign{\smallskip}\hline
	\end{tabular}}
\end{table*}

\subsection{Discussion}
The posterior computed with the 2D subspace has changed greatly from the prior in the dominant parameters in the first two eigenvectors $\eigVec_1$ and $\eigVec_2$.
This can be observed by large correlation between the first and second and the sixth and eighth parameters.
Since the subsequent three eigenvectors $\eigVec_3$, $\eigVec_4$, and $\eigVec_5$ have potentially significant eigenvalues, inference in the 5D subspace leads to additional significant changes.
Using the 5D subspace better infers more directions of the parameter space (and individual parameters), which could be valuable in certain cases (e.g. predicting a different quantity of interest).
In particular, parameters four, five, and eight are much better inferred by the 5D space.
This can be explained by again regarding the active variable $\y$ as a linear combination.
Taking the first eigenvectors $\eigVec_1,\ldots,\eigVec_5$ into account for $\y$ leads to putting weights on more components of $\x$ than in the 2D case, i.e. more components and their multivariate behavior are informed.

\begin{figure}
	\centering
	\includegraphics[scale=0.5]{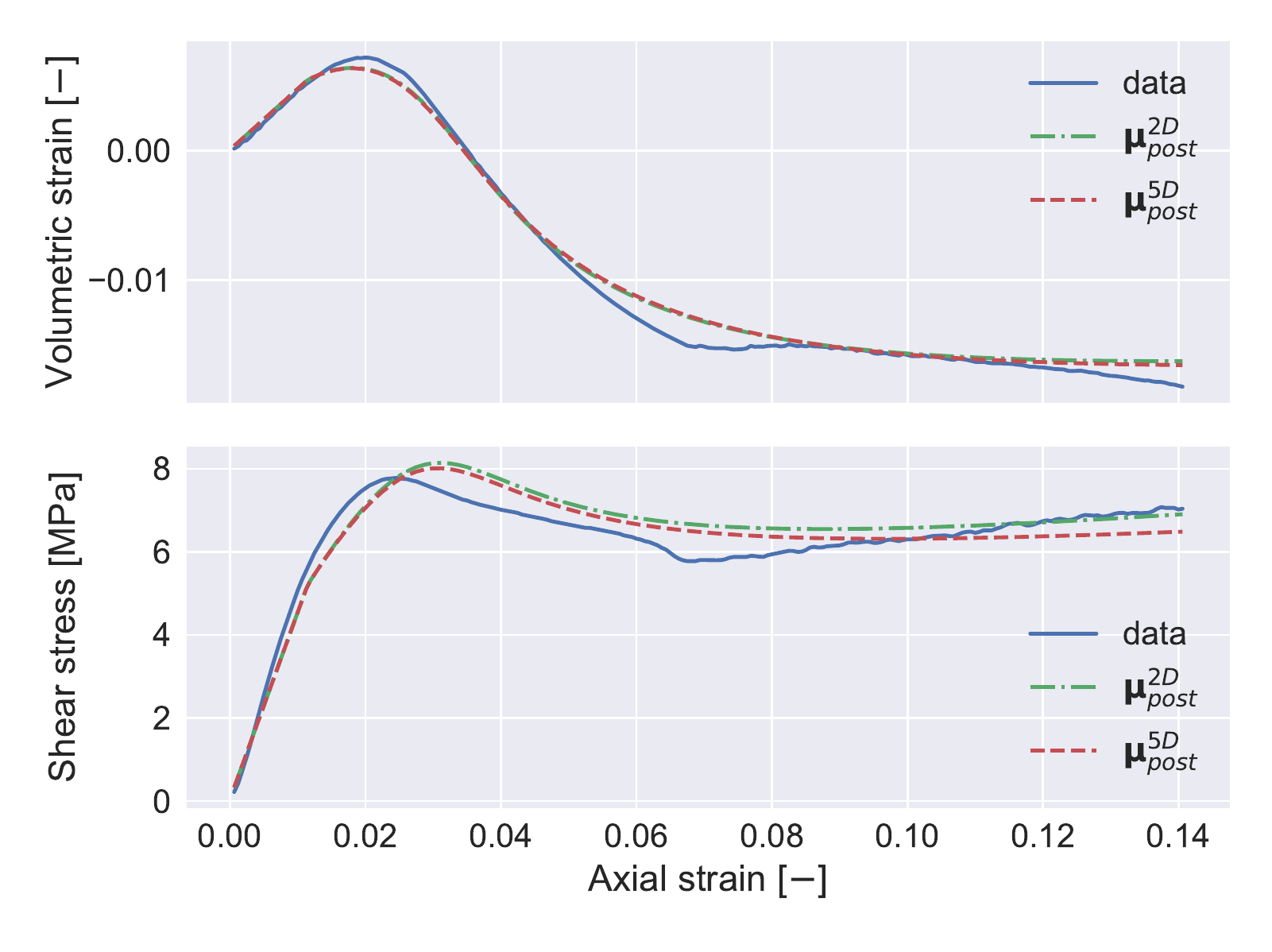}
	\caption{Model evaluations for posterior means computed with the 2D and 5D subspace.}
	\label{fig:posterior_mean}
\end{figure}

We additionally evaluate the forward model with the posterior means as inputs to recover mean stress-strain curves (see Fig.~\ref{fig:posterior_mean}).
Both model evaluations have similar curves, which match the experimental data relatively well.
While the 2D and 5D active subspaces give different posteriors, the data misfit is not very sensitive in directions in which they differ.
The data misfit varies much more in the directions of the first two eigenvectors which define the 2D active subspace and also are components of the 5D active subspace.
Hence, the posteriors from the 2D and 5D subspaces both are reasonable for the given quantities of interest.
In future work, it may be possible to construct alternative Bayesian inference problems (with different quantities of interest, noise levels, and weights) that result in posteriors that better match important physical characteristics of the shear stress and volumetric strains such as peak locations, peak heights, and second hardening rates and locations.

\section{Conclusions}
\label{sec:conclusions}

In this study, we present a state-of-the-art constitutive model describing the geomechanical behavior of methane hydrate bearing sands. 
The model builds on the Drucker-Prager theory for soil plasticity and aims to capture the distinct secondary hardening phase, observed in our drained triaxial experiments, through strain-dependent evolution laws for the friction and dilation parameters.
Little has been reported in literature on the secondary hardening in laboratory testing of gas hydrates, and as such, our current understanding of this effect and its implications on the field-scale geomechanical behavior of gas hydrate geosystems is very limited.
This lack of knowledge is directly reflected in our poor intuition for the model parameters and makes the tasks of parameter estimation and model calibration particularly challenging.
In this work, we focused on a single experimental data-set, essentially a ``training" set, to test the hypothesis that there exist lower-dimensional structures (i.e. active subspaces) in the space of parameters of our methane hydrate constitutive model which can be utilized to efficiently perform parameter estimation in a Bayesian setting.
The bases of the active subspaces consist of linear combinations of model parameters, providing insights into the relationships between these parameters and how they affect the model output. 
Such insights are extremely helpful in identifying the dominant parameters (or hyper-parameters).

Out of the full 8D space of uncertain parameters in our model, we are able to identify 2D and 5D active subspaces in which the model inferences are most sensitive.
The 2D and 5D active subspaces are used to develop efficient and accurate regression-based surrogate models.
Because of the high computational cost of the methane hydrate model and the non-trivial dimension of the full space of uncertain parameters, traditional MCMC methods are computationally prohibitive with the full model, but the surrogates allowed for an efficient MCMC algorithm for Bayesian inference.
The mean of the posterior density calculated from the algorithm matched extremely well with the observed experimental data.
Furthermore, in the 2D active subspace, the cohesion parameter $c$ and the initial residual friction parameter $\alpharesmin$ emerge as dominant, which matches very well with the expectation from any coulomb-type plasticity model where the initial yield surface (i.e. $q+\alpharesmin-c=0$) controls the onset of plasticity.
In general, the parameters $c$ and $\alpharesmin$ can be estimated with relatively high confidence through triaxial testing, and in future investigations our focus will be on identifying active subspaces for the other parameters for which we indeed lack direct estimation procedures.
Additionally, global sensitivity metrics for parameters, which are also constructible with active subspaces \cite{constantine2017global}, are of interest as well as the inference of parameters for other types of sediment saturations.

Active subspace analysis is a relatively new tool which has only recently been applied to stochastic inversion.
If active subspaces exist in a space of uncertain parameters, they can be used to effectively reduce the dimension of parameter spaces and to accelerate Bayesian inference.
Active subspace analysis should be considered as a potential tool in a wide range of geophysical applications, because they may be utilized for highly efficient parameter estimation which would otherwise be computationally prohibitive.

\section*{Acknowledgements}
Financial support for BW, SM, and MTP  was provided by the German Research Foundation (DFG, Project WO 671/11-1).
The work of SG and CD was further funded by the German Federal Ministries of Economy (BMWi) and Education and Research (BMBF) through the SUGAR project (grant no. 03SX250, 03SX320A \& 03G0856A), 
and the EU-FP7 project MIDAS (grant agreement no. 603418).

\end{document}